\documentclass[reqno, 12pt]{amsart}
\usepackage{amsfonts}
\usepackage{amssymb,amsmath}
\usepackage{amsthm}
\usepackage{upref}

\makeatletter
\def\LaTeX{\leavevmode L\raise.42ex
    \hbox{\kern-.3em\size{\sf@size}{0pt}\selectfont A}\kern-.15em\TeX}
 \makeatother

\numberwithin{equation}{section}

\newtheorem{lemma}{Lemma}[section]
\newtheorem{theorem}[lemma]{Theorem} 
\newtheorem{corollary}[lemma]{Corollary}
\newtheorem{proposition}[lemma]{Proposition}

\theoremstyle{definition}

\newtheorem{definition}[lemma]{Definition}
\newtheorem{example}[lemma]{Example}

\newtheorem{remark}[lemma]{Remark}

\newcommand{\Ker}{\operatorname{Ker}}

  \newcommand{\e}{\eqref}

\newcommand{\q}{\quad}

\newcommand{\ti}{\tilde}
\newcommand{\wt}{\widetilde}
\newcommand{\wh}{\widehat}
\newcommand{\la}{\langle}
\newcommand{\ra}{\rangle}

\newcommand{\ov}{\overline}
 \renewcommand{\d}{\delta}

  \newcommand{\rank}{\operatorname{rank}}

\renewcommand\Im{\operatorname{Im}}
\renewcommand\Re{\operatorname{Re}}

\newenvironment{pf}{\begin{proof}}{\end{proof}}

\def\qqq{\mathrel{\subset\mkern-15mu\lower.38ex\hbox{${\scriptscriptstyle\rightarrow}$}}}
\let\cal\mathcal

\let\Bbb\mathbb

  \DeclareMathOperator{\spec}{spec}

\begin{document}

\title{Criteria for  Hankel  operators to be sign-definite}
\author{ D. R. Yafaev}
\address{ IRMAR, Universit\'{e} de Rennes I\\ Campus de
  Beaulieu, 35042 Rennes Cedex, FRANCE}
\email{yafaev@univ-rennes1.fr}
\keywords{Hankel  operators, convolutions,   necessary and sufficient conditions for    the  positivity,  the sign-function,  operators of finite rank, the Carleman operator and its perturbations}
\subjclass[2000]{47A40, 47B25}

%  \date{\today}

\begin{abstract}
We show that total multiplicities of negative and positive spectra of a self-adjoint Hankel  operator $H$ with kernel $h(t)$ and of an operator of multiplication by some real function $s(x)$
  coincide.  In particular,  $\pm H\geq 0$ if and only if $\pm s(x)\geq 0$. The kernel $h(t)$ and its ``sign-function" $s(x)$ are related by an explicit formula.  An expression of $h(t)$ in terms of $s(x)$ leads to an exponential representation of $h(t)$.
  Our approach directly applies to various classes of   Hankel  operators. In particular, for Hankel  operators of  finite rank, we  find an explicit formula for the total multiplicity of their negative and positive spectra. 
    \end{abstract}

\maketitle

% \thispagestyle{empty}

%************************************************************
\section{Introduction}  
%***********************************************************

{\bf 1.1.}
Hankel operators can be defined as integral operators
\begin{equation}
(H f)(t) = \int_{0}^\infty h(t+s) f(s)ds 
\label{eq:H1}\end{equation}
in the space $L^2 ({\Bbb R}_{+}) $ with kernels $h$ that depend  on the sum of variables only. Of course $H$ is symmetric if $  h(t)=\ov{h(t)}$. There are very few cases where Hankel operators can be explicitly diagonalized. We mention classical results by F.~Mehler \cite{Me}, T.~Carleman \cite{Ca}, W.~Magnus \cite{Ma} and M.~Rosenblum \cite{Ro}. They are treated in a unified way in \cite{Y1} where some new examples are also considered. 

Our goal here is to find   necessary and sufficient conditions for the positivity\footnote{We   always use the term ``positive" instead of a more precise but lengthy term  ``non-negative."} of Hankel operators.  This question seems to be of particular importance because of its intimate relation to a representation of the function $h(t)$ as the Laplace transform of a (positive) measure. Such representations are continuous analogues of
moment problems.  According to the Hamburger theorem (see, e.g., the book \cite{AKH}) the positivity of a discrete Hankel operator is equivalent to the existence of a solution of the corresponding moment problem.

 \medskip
 
 {\bf 1.2.} 
 Our condition of the positivity\footnote{We   usually discuss conditions for $H\geq 0$, but of course replacing $H$ by $-H$ we obtain conditions for $H\leq 0$.} of Hankel operators   is quite explicit.  
Let    $B$,
 \begin{equation}
(B g) (\xi) =  \int_{- \infty}^\infty   b(\xi-\eta) g(\eta)d\eta,
\label{eq:BB}\end{equation}
 be the operator in the space $L^2 ({\Bbb R})$ of the convolution with the function 
  \begin{equation}
 b(\xi) = \frac{1}{ 2\pi}  \frac{ \int_{0}^\infty   h (t) t^{ -i \xi}   dt}{ \int_{0}^\infty   e^{-t} t^{ -i \xi}   dt} .
\label{eq:Bb}\end{equation}
Of course $  b(-\xi)=\ov{ b(\xi)}$ if $  h(t)=\ov{ h(t)}$ so that the operator $B$ is symmetric.
Our main result is that the Hankel operator $H\geq 0$ if and only if $B\geq 0$. We call $b(\xi)$ the $b$-{\it function} of the Hankel operator $H$ (or of the kernel $h(t)$), and we use the term {\it  the sign-function} or $s$-{\it function}  for its  Fourier transform $s(x)$. So our result, roughly speaking,  means that a Hankel operator $H$ is positive if and only if its   sign-function $s(x)$ is positive.  Note that in specific examples we consider, functions $s(x)$ may be of a quite different nature. For instance, $s(x)$ may be a polynomial or, on the  contrary, it may be a distribution; for example, it may be a combination of delta functions and their derivatives.

Our proofs rely on the identity
  \begin{equation}
 (H f,f)=(B g,g) 
\label{eq:MAID}\end{equation}
where
 \begin{equation}
g(\xi)=\Gamma(1/2+i \xi)\ti{f}(\xi)=:( \Xi f) (\xi), 
\label{eq:MAID1}\end{equation} 
  $\ti{f}(\xi)$ is the Mellin transform of $f (t)$ and $\Gamma(\cdot)$ is the gamma function.     Actually,  the identity  \e{eq:MAID} or, equivalently,
   \begin{equation}
 H= \Xi^* B \Xi
\label{eq:MAIDs}\end{equation}
allows us to find the numbers of positive and negative eigenvalues
  of a Hankel operator $H$. For a self-adjoint operator $A$, we denote by $N_{+}(A)$ (by $N_{-}(A)$) the total mutiplicity of its strictly positive (negative) spectrum. Then  the identity  \e{eq:MAIDs}
shows that 
 \begin{equation}
N_{\pm}(H)= N_{\pm}(B). 
\label{eq:MAIDs1}\end{equation}

This result can be compared with Sylvester's inertia theorem  which states the same   for Hermitian matrices $H$ and $B$ related by equation \e{eq:MAIDs} provided the matrix $\Xi$ is invertible. In contrast to the linear algebra, in our case the operators $H$ and $B$ are of a completely different nature and $B$ (but not $H$) admits an explicit spectral analysis.

Of course the results above can be reformulated in terms of the sign-function $s(x)=\sqrt{2\pi} (\Phi^* b)(x)$ where $\Phi$ is the Fourier transform. Let us   introduce the operator $S=\Phi^* B\Phi$ of multiplication by $s(x)$. Then \e{eq:MAIDs} means that
$ H= \wh{\Xi}^* S \wh{\Xi}$ where $\wh{\Xi}=\Phi^*\Xi$, and \e{eq:MAIDs1} means that
$N_{\pm}(H)= N_{\pm}(S)$.
  
  \medskip
 
 {\bf 1.3.} 
 The precise meaning of formula \e{eq:Bb} requires some discussion.  Observe  that the denominator in \e{eq:Bb} coincides with the numerator for the special case $h(t)=e^{-t}$.
 It equals $\Gamma(1- i \xi)$ and hence exponentially tends to zero as $|\xi|\to\infty$. Therefore $b(\xi)$ is a ``nice" function of $\xi$ only under very restrictive assumptions on the kernel $h(t)$. Thus to cover natural examples, we have, on the contrary, to extend a class of kernels   and to work with distributions $h(t)$. The choice of appropriate spaces of distributions is also very important. In order to be able to divide in \e{eq:Bb} by an exponentially decaying function,   we assume that the numerator belongs to the class of distributions $C_{0}^\infty ({\Bbb R})'$. It means that the Fourier transform of    the function    
$ \theta (x)= e^x h(e^x)$ should belong  to $C_{0}^\infty ({\Bbb R})'$. The corresponding class of functions $h(t)$ will be denoted ${\cal Z}_{+}'$.  Under the assumption $h\in {\cal Z}_{+}'$ we   have $b \in C_{0}^\infty ({\Bbb R})'$. The Schwartz space ${\cal S}({\Bbb R})'$ is too restrictive for our purposes which is seen already on the example of finite rank Hankel operators. 

The condition $ \Phi \theta \in C_{0}^\infty ({\Bbb R})'$ for the validity of 
 the identity  \e{eq:MAID} is very general. It is satisfied for {\it all} bounded, but also for a wide class of unbounded, Hankel operators $H$. More than that, it is not even required that $H$ be defined by formula \e{eq:H1}     on some dense set. Therefore we work with quadratic forms $(Hf,f)$ which is more convenient and yields more general results.  In this context it is   natural to consider distributions $h(t)$ which makes the theory   self-consistent.

   If however $h\in L^1_{\rm loc}({\Bbb R}_{+})$, then $ \theta\in {\cal S}({\Bbb R})'  $ and hence $\Phi \theta\in    C_{0}^\infty ({\Bbb R})'$ if
 \begin{equation}
\int_{0}^\infty | h(t)| (1+| \ln t |)^{-\kappa } dt< \infty
 \label{eq:ass}\end{equation}
 for some $\kappa $.  Condition \e{eq:ass} is quite general; it also does not  require that the corresponding Hankel operator be bounded. For example, it admits kernels 
  \begin{equation}
h(t)= P (\ln t)t^{-1}
 \label{eq:LOG}\end{equation}
where  $P (x)$
is an arbitrary polynomial. Note that Hankel operators with such kernels are   bounded for $P (x)={\rm const}$ only.

As far as test functions $f(t)$ are concerned, we require that their Mellin transforms $\ti{f} \in C_{0}^\infty ({\Bbb R})$. Then both sides of  \e{eq:MAID} are well defined  and    the identity    holds.   
We note that the distribution  $b(\xi)$ is ``worse" than the kernel $h (t)$. On the contrary, due to the factor $\Gamma(1/2+i \xi)$ in \e{eq:MAID1}, the test function $g(\xi)$
is ``better" than $f(t)$.  In the case of bounded operators $H$, this permits us to extend  the main  identity  \e{eq:MAID} to all elements $f\in L^2 ({\Bbb R}_{+})$.

 \medskip
 
 {\bf 1.4.} 
 It turns out that the knowledge of the sign-function $s(x)$ allows one to recover the kernel $h(t)$ by the formula
  \begin{equation}
  h (t)  =    \int_{0}^\infty     e^{-t \lambda} h^\natural (\lambda)  \lambda d \lambda
   \label{eq:conv1}\end{equation}
   where
  \begin{equation}
 h^\natural (\lambda)=  \lambda^{-1} s (-\ln \lambda).
 \label{eq:B}\end{equation} 
Thus     $h (t)$ is the Laplace transform of the distribution $ \lambda h^\natural  (\lambda)$. It is noteworthy that $h^\natural \in{\cal Z}_{+}'$ and the correspondence $h\mapsto h^\natural$ is a continuous one-to-one mapping of ${\cal Z}_{+}'$ onto itself.

  Representation \e{eq:conv1} does not require the positivity of $H$. If however $H\geq 0$, then combining our results with the Bochner-Schwartz theorem, we obtain that $\lambda h^\natural (\lambda)  d \lambda=   dm(\lambda)$ where $dm(\lambda)$ is a positive measure. In this case 
 \begin{equation}
  h (t)  =    \int_{0}^\infty     e^{-t \lambda}    d m(\lambda).
   \label{eq:Conv}\end{equation}
   This representation implies that $h(t)$ is necessarily a completely monotonic function and, in particular, $h\in C^\infty ({\Bbb R}_{+})$. Note that the converse statement is also true: a completely monotonic function admits representation \e{eq:Conv}. This is one of famous Bernstein's theorems (see his paper \cite{Bern}, or the book by  N.~Akhiezer \cite{AKH} or the book by  D.~V.~Widder \cite{Wid}). In contrast to the
   Bernstein  theorem, we deduce representation \e{eq:Conv} from the positivity of the Hankel operator with kernel $h(t)$ and show that the measure  $dm(\lambda)$ satisfies  for some $\varkappa$ the condition
         \begin{equation}
\int_{0}^\infty (1+|\ln\lambda|)^{-\varkappa } \lambda^{-1} dm (\lambda )<\infty.
\label{eq:Sch1}\end{equation}

 Hankel operators can   also be realized as operators in the space   of sequences $l_{+}^2$. The relation between this discrete representation and the continuous  representation we consider is given by the unitary transformation of  $ l_{+}^2$ onto $ L^2 ({\Bbb R}_{+})$ constructed in terms of the Laguerre polynomials. Thus all our results can in principle be translated into the discrete representation. This hopefully will be discussed in another article on this subject.

  \medskip
 
 {\bf 1.5.} 
 A large part of the paper is devoted to applying the general theory to various classes  of Hankel operators although we do not try to cover all possible cases. 
 In some examples the   sign definiteness of $H$ can   also be verified or refuted with the help of Bernstein's theorems. Note however that our approach yields additionally an explicit formula for the total numbers of  negative and positive eigenvalues of $H$. 
  
In Section~5,  we present  such a   formula for  Hankel operators  of  finite rank.  
 Then we consider   two specific examples.
  The first one is given by the formula
 \begin{equation}
h(t)= t^k e^{-\alpha t}, \q \alpha> 0, \, k \geq -1.
\label{eq:E1r}\end{equation}
  Note that the Hankel operator $H$ with such kernel has finite rank for   $k\in{\Bbb Z}_{+}$ only.
We show that    $H$   is positive if and only if $k\leq 0$.
The second class of kernels is defined by the formula  
   \begin{equation}
 h(t)=e^{-t^r},\q r>0.
\label{eq:F5}\end{equation}
It turns out that the corresponding Hankel  operator is positive if and only if $r\leq 1$.  

Section~6 is devoted to a study of Hankel  operators $H$ with non-smooth kernels. In this case both numbers $N_{\pm} (H)$ are infinite, and we find the asymptotics of eigenvalues of $H$.

Finally, in Section~7 we consider perturbations of the Carleman operator, that is, of
the Hankel operator with   kernel $h_{0}(t)=  t^{-1}$ by various classes of compact Hankel operators. The Carleman operator can be explicitly diagonalized by the Mellin transform.
We recall that it has the absolutely continuous spectrum $[0,\pi]$ of multiplicity $2$. The Carleman operator plays the distinguished role in the theory of Hankel operators. In particular, it is important for us   that its sign-function  $s(x)=1$.
As was pointed out by J.~S.~Howland in \cite{Howland},   Hankel  operators are to a certain extent similar to differential operators. 
In terms of this analogy, the Carleman operator $H_{0}$   plays the role of the ``free" Schr\"odinger operator $D^2$, $D=-i d /dx $,  in the space $L^2 ({\Bbb R}) $. Furthermore, Hankel operators $H$ with ``perturbed" kernels $h (t)=t^{-1}+v(t)$  can be compared to Schr\"odinger operators  $D^2+ {\sf V}(x)$. The assumption that
   $v(t)$    decays sufficiently rapidly as $t\to \infty$ and is not too singular as $t\to 0$ corresponds to a sufficiently rapid decay 
 of  the potential $  {\sf V}(x)$  as $|x|\to \infty$.   
 
As shown in \cite{Y2},   the results     on the discrete spectrum of the operator $H$ lying {\it above} its essential spectrum $[0,\pi]$ are close in spirit to the results on the discrete spectrum of the Schr\"odinger operator  $D^2+ {\sf V}(x)$. On the contrary,  the results on the negative spectrum of the Hankel operator $H $  are drastically different.
   In particular,  contrary  to the case of differential operators with decaying coefficients, the finiteness of the negative spectrum of the Hankel operator $H $  is not determined by the behaviour of $v(t)$ at singular points $t=0$ and $t=\infty$. As an example, consider the   Hankel operator  with kernel 
      \[
      h(t)=t^{-1}-\gamma e^{-t^r}, \q r \in (0,1).
      \]
    Now the kernel of the perturbation is   the function which decays faster than any power of $t^{-1}$ as $t\to \infty$ and it has the finite limit as $t\to 0$. Nevertheless we show that the negative spectrum   of $H$ is infinite if $\gamma> \gamma_{0}$ (here $\gamma_{0}=\gamma_{0} (r)$ is an explicit constant) while $H$ is positive if $\gamma\leq \gamma_{0}$. Such a phenomenon has no analogy for   Schr\"odinger operators with   decaying potentials. However it occurs for   three-particle Schr\"odinger operators and is known as the Efimov effect.

   We also study   perturbations of the Carleman operator $H_{0}$ by Hankel operators $V$ of  finite rank. Here we obtain a striking result: the total numbers of   negative   eigenvalues of the operators $H=H_{0}+V$ and $V$ coincide.

 As examples, we consider   only bounded Hankel operators in  this paper. However, our general results directly apply to a wide class of unbounded   operators, such as Hankel operators with kernels \e{eq:LOG}. Moreover, with slight modifications our method works also for kernels   \e{eq:E1r} where $\alpha\geq 0$ and $k$ is an {\it arbitrary} negative number. In this case condition  \e{eq:ass} is not satisfied. Hankel operators with kernels \e{eq:LOG} and \e{eq:E1r} generalize the Carleman operator, and we call them quasi-Carleman operators. They will be studied elsewhere.

   \medskip
 
 {\bf 1.6.}
 Let us briefly describe the structure of the paper. We obtain the main identity \e{eq:MAID} and the reconstruction formula \e{eq:conv1} in Section~2. 
 Necessary information on bounded Hankel operators (including a continuous version of the Nehari theorem) is collected  in Section~3. 
 In Sections~2 and 3 we   do not assume that the function $h$ is real, i.e., the corresponding Hankel operator $H$ is not necessarily symmetric. 
  Spectral consequences of the formula \e{eq:MAID} and, in particular, criteria of the sign-definiteness of Hankel operators are formulated in Section~4. 
  In Sections~5, 6 and 7   we apply the general theory to particular classes of Hankel operators.

    Let us introduce some standard
 notation. We denote by $\Phi$, 
\[
(\Phi u) (\xi)=  (2\pi)^{-1/2} \int_{-\infty}^\infty u(x) e^{ -i x \xi} dx,
\]
  the Fourier transform.       The space $\cal Z= \cal Z ({\Bbb R})$ of test functions is defined  as the subset  
 of the Schwartz  space ${\cal S}={\cal S} ({\Bbb R}) $ which consists of functions $\varphi  $ admitting the analytic continuation to   entire functions in the   complex plane $\Bbb C$   and satisfying    bounds 
 \[
  | \varphi (z)| \leq C_{n}  (1+| z |)^{-n} e^{r |\Im z |},\q \forall z\in \Bbb C,
  \]
  for some $r=r(\varphi)>0$  and all $n$. We recall that the Fourier transform $\Phi: C_{0}^\infty ({\Bbb R}) \to \cal Z$ and $\Phi^* : \cal Z \to C_{0}^\infty ({\Bbb R})$.  The dual classes of distributions (continuous antilinear functionals) are denoted ${\cal S}'$, $C_{0}^\infty ({\Bbb R})'$ and ${\cal Z}'$, respectively.
  In general, for a linear topological space $\cal L$, we use the notation ${\cal L}'$ for its dual space. The Dirac function is standardly denoted $\d(\cdot)$.

  We use the notation   ${\pmb\la} \cdot, \cdot {\pmb\ra}$ and $\la \cdot, \cdot\ra$ for scalar products and duality symbols in $L^2 ({\Bbb R}_{+})$ and $L^2 ({\Bbb R})$, respectively. They are always linear in the first argument and antilinear in the second argument. The letter $C$ (sometimes with indices) denotes various positive constants whose precise values are inessential.

%************************************************************
\section{The main identity  }  
%*

{\bf 2.1.}
Let us consider a Hankel operator $H$ defined by equality \e{eq:H1} in the space $L^2 ({\Bbb R}_{+})$.       Actually, it is more convenient to work with sesquilinear forms $(Hf_{1} ,f_{2})$ instead of operators. 
 
 Before giving precise definitions, let us explain our construction at  a formal level.
   It follows from \e{eq:H1} that
\begin{align}
(Hf_{1} ,f_{2})= &\int_{0}^\infty \int_{0}^\infty h(t+s) f_{1}(s)\overline{f_{2}(t)} dtds
\nonumber\\
= & \int_{0}^\infty  h(t ) \overline{F(t)} dt=: {\pmb\la} h,  F {\pmb\ra},
 \label{eq:HH}\end{align}
 where  
\begin{equation}
F(t)=   \int_{0}^t    \overline{f_{1}(s)} f_{2}(t-s) ds =: ( \bar{f}_{1}\star f_{2})(t)
 \label{eq:HH1}\end{equation}
 is the   Laplace convolution of the functions $ \bar{f}_{1}$ and $ f_{2}$.   Formula  \e{eq:HH} allows us to consider  $h$   as a distribution with the test function $F$ defined by  \e{eq:HH1}.   Thus the Hankel form will be defined by the relation
  \begin{equation}
h[f_{1},f_{2}]= {\pmb\la} h, \bar{f}_{1}\star f_{2} {\pmb\ra}.
 \label{eq:HH2}\end{equation}

 Let us introduce the test function
 \begin{equation}
\Omega(x)=F(e^x) =: ({\cal R} F)(x)
 \label{eq:HH3s}\end{equation} 
 and the distribution  
 \begin{equation}
 \theta (x)= e^x h(e^x)
 \label{eq:HH3}\end{equation}
 defined for $x\in{\Bbb R}$.
  Setting in  \e{eq:HH} $t=e^x$, we see that
 \begin{equation}
  {\pmb\la} h, F {\pmb\ra} = \int_{-\infty}^\infty \theta (x ) \overline{\Omega(x)} d x=: { \la} \theta, \Omega { \ra}.
 \label{eq:HH4}\end{equation}

 We are going  to consider the form \e{eq:HH4} on pairs $F, h$ such that the corresponding test function $\Omega$ defined by \e{eq:HH3s}  is an element of the space $\cal Z$  of analytic functions and the corresponding distribution $\theta$ defined by  \e{eq:HH3} is an element of the dual space $\cal Z'$. The set of all such $F$ and $h$ will be denoted ${\cal Z}_{+}$ and ${\cal Z}_{+}'$, respectively, that is,
 \begin{equation}
 F\in {\cal Z}_{+} \Longleftrightarrow \Omega\in  {\cal Z}
 \q {\rm and}\q
h\in {\cal Z}_{+}' \Longleftrightarrow \theta\in  {\cal Z}'.
\label{eq:heta}\end{equation}
Of course, the topology in ${\cal Z}_{+}$ is induced by that in $\cal Z$ and ${\cal Z}_{+}'$ is dual to ${\cal Z}_{+}$. Note that $h\in{\cal Z}_{+}'$ if $h\in L^1_{\rm loc}({\Bbb R}_{+})$ and integral
 \e{eq:ass} is convergent
 for some $\kappa $. In this case the corresponding function \e{eq:HH3} satisfies the condition
  \[
\int_{-\infty}^\infty |\theta (x)| (1+| x |)^{-\kappa  } dx< \infty,
 \]
 and hence $\theta \in{\cal S}'\subset {\cal Z}'$.
  
Define  the unitary operator $U: L^2 ({\Bbb R}_{+})\to L^2 ({\Bbb R} )$  by the equality
\begin{equation}
(U f)(x) =e^{x/2} f(e^x).
\label{eq:M1X}\end{equation}
Let the set $\cal D$ consist of functions $f  (t)$  such that $  Uf  \in {\cal Z}$. Since
\[
f(t)= t^{-1/2}  (Uf) (\ln t)
\]
and ${\cal Z}\subset {\cal S}$, we see that   functions $f  \in \cal D$ and their derivatives satisfy the estimates
 \[
|f^{(m)}(t)|=  C_{n,m} t^{-1/2-m} (1+|\ln t|)^{-n}
\]
for all $n$ and $m$. Obviously, $f \in \cal D$ if and only if $\varphi(t)=t^{1/2} f(t)$ belongs to the class ${\cal Z}_{+} $.

 Let us show that
  form \e{eq:HH2} is correctly defined on     functions $f_{1}, f_{2}\in \cal D$. To that end, we have to verify that     function \e{eq:HH1}  belongs to the space $ {\cal Z}_{+}$ or, equivalently,   function \e{eq:HH3s} belongs to the space $\cal Z$. This requires some preliminary study which will   also allow us to derive a convenient representation for   form \e{eq:HH2}.

   Recall that the Mellin transform ${\bf M}: L^2 ({\Bbb R}_{+})\to L^2 ({\Bbb R} )$ is defined by the formula
\begin{equation}
({\bf M} f) (\xi)=  (2\pi)^{-1/2} \int_{0}^\infty f(t) t^{-1/2-i\xi} dt.
\label{eq:M1}\end{equation}
Of course, ${\bf M}=\Phi U$ where $\Phi$ is the Fourier transform 
 and $U $ is   operator \e{eq:M1X}.
Since both $\Phi$ and $U$ are unitary, the operator ${\bf M}$ is also unitary.
The inversion of the formula \e{eq:M1} is given by the relation 
\begin{equation}
f (t)=  (2\pi)^{-1/2} \int_{-\infty}^\infty \tilde{f}(\xi) t^{-1/2+i\xi} d\xi, \q \tilde{f}= {\bf M}f.
\label{eq:M2}\end{equation}

 Let $\Gamma(z)$ be the gamma function. Recall that $\Gamma(z)$ is a holomorphic function in the right half-plane and  $\Gamma(z)\neq 0$ for all $z\in{\Bbb C}$.     According to  the Stirling formula  the   function $\Gamma(z)$ tends to zero exponentially as $|z|\to\infty$ parallel with the imaginary axis. To be more precise,
    we have
\begin{equation}
   \Gamma(\alpha +i \lambda) = e^{\pi i (2\alpha-1)/4} (2\pi / e)^{1/2}  \lambda^{\alpha-1/2} e^{i\lambda (\ln \lambda -1)} e^{-\pi \lambda /2}\big(1+O(\lambda^{-1})\big)  
\label{eq:M11}\end{equation}
for a fixed $\alpha>0$ and $ \lambda\to +\infty$.
Since $ \Gamma(\alpha - i \lambda)=\overline{\Gamma(\alpha +i \lambda)}$, this yields also the asymptotics of $   \Gamma(\alpha +i \lambda)$ as $\lambda\to -\infty$.

 If $f_{j} \in \cal D$, $j=1,2$,  then   $\ti{f}_{j}= {\bf M} f_{j}= \Phi U f_{j}\in C_{0}^\infty ({\Bbb R})$ and hence the functions
\begin{equation}
 g_{j}(\xi) =    \Gamma(1/2+i \xi) \tilde{f}_{j}(\xi), \q j=1,2,
\label{eq:M10}\end{equation}  
    also belong to the class $ C_{0}^\infty ({\Bbb R})$.
     Let us introduce the convolution 
      of the functions $g_1 $ and $g_2 $,
 \[
 (g_1 * g_{2})(\xi) = \int_{-\infty}^\infty  g_{1}( \xi- \eta)  g_2(\eta ) d\eta,
\]
and set
 \[
 ({\cal J} g)(\xi)=g(-\xi).
 \]
 We have the following result.

    \begin{lemma}\label{fg1}
 Suppose   that $f_{j}  \in \cal D$, $j=1,2$,  and define   functions $g_{j}(\xi)$ by equality \e{eq:M10}.  Let the function $\Omega(x)$ be defined by formulas \e{eq:HH1} and \e{eq:HH3s}.  Then  
 \begin{equation}
(\Phi \Omega)(\xi)= (2\pi)^{-1/2  } \Gamma (1+ i\xi)^{-1}  (({\cal J}\bar{g}_1) * g_{2}) (\xi) .
\label{eq:M81}\end{equation}
 \end{lemma}
 
 \begin{pf}
  Substituting \e{eq:M2} into definition \e{eq:HH1}, we see that
\[
F(t)= (2\pi)^{-1 }\int_{0}^t ds 
\int_{-\infty}^\infty \overline{\tilde{f}_1(\tau)} (t-s)^{-1/2 - i \tau}d\tau
\int_{-\infty}^\infty \tilde{f}_2 (\sigma) s^{-1/2+ i\sigma}d \sigma.
\]
Observe that
\[
\int_{0}^t       (t-s)^{-1/2 - i \tau}  s^{-1/2+ i \sigma}  ds= t^{  i( \sigma-\tau)}  \frac{\Gamma(1/2-i\tau)
 \Gamma(1/2+ i \sigma)}{ \Gamma (1+ i(\sigma-\tau))}.
\]
Then using definition  \e{eq:M10} 
we obtain the representation
 \begin{align*}
F(t)= & (2\pi)^{-1 }  \int_{-\infty}^\infty  \int_{-\infty}^\infty  t^{  i( \sigma-\tau)} \Gamma (1+ i(\sigma-\tau))^{-1} \overline{ {g}_{1}(\tau)}   {g}_{2}(\sigma)     d\tau d \sigma
\\
= & (2\pi)^{-1  }  \int_{-\infty}^\infty   t^{  i \xi  } \Gamma (1+ i\xi)^{-1} ( ({\cal J}\bar{g}_1) * g_{2})(\xi)  d\xi,
 \end{align*}
 whence
\[
\Omega(x)= (2\pi)^{-1  }  \int_{-\infty}^\infty   e^{  i \xi x} \Gamma (1+ i\xi)^{-1} (({\cal J}\bar{g}_1) * g_{2}) (\xi)  d\xi.
\]
This  
  is equivalent to formula \e{eq:M81}. 
    \end{pf}

Observe that the function $\Gamma (1+ i\xi)^{-1}$ in the right-hand side of \e{eq:M81} tends to infinity exponentially as $|\xi|\to\infty$. Nevertheless $\Phi \Omega\in C_{0}^\infty ({\Bbb R})$ because $({\cal J}\bar{g}_1) * g_{2} \in C_{0}^\infty ({\Bbb R})$  for $ g_1, g_{2}\in C_{0}^\infty ({\Bbb R})$. Thus we have
    
        \begin{corollary}\label{fg1c}
 Let  $f_{j} \in \cal D$, $j=1,2$, and let the function $\Omega(x)$ be defined by formulas \e{eq:HH1} and \e{eq:HH3s}. Then $\Omega\in\cal Z$ or, equivalently, $F \in {\cal Z}_{+}$.  
 \end{corollary}  
    
Now we are in a position to give the precise definition.

  \begin{definition}\label{HH}
Let   $h\in  {\cal Z}_{+}'$ and $f_{j} \in \cal D$, $j=1,2$.   Then the Hankel sesquilinear form is defined by the relation \e{eq:HH2}. 
 \end{definition}
 
 We shall see in subs.~2.4 that $h\in  {\cal Z}_{+}'$ is determined uniquely by the values $ { \pmb\la} h ,  \bar{f}_{1}\star f_{2} {\pmb \ra}$ on $f_1, f_2 \in \cal D$, that is, $h=0$  if $ { \pmb\la} h ,  \bar{f}_{1}\star f_{2} {\pmb \ra}=0$ for all $f_1, f_{2} \in \cal D$.
 
 Of course   definition \e{eq:HH2} can be rewritten as
    \begin{equation}
h[f_{1}, f_{2}]  = { \la}\theta , \Omega { \ra} 
 \label{eq:hhfx}\end{equation}
 where
 \[
   \Omega (x)= (\bar{f}_{1}\star f_{2})(e^x)
\]
    and $\theta$  is distribution     \e{eq:HH3}.
 
 We sometimes write $h[f_{1},f_{2}]$ as integral  \e{eq:HH} keeping in mind that its precise meaning is given by Definition~\ref{HH}.

    \medskip
 
 {\bf 2.2.}
 Our next goal is  to show that \e{eq:hhfx} is the sesquilinear form of the convolution operator $B$, that is, it equals the right-hand side of  \e{eq:MAID}. Here the representation of Lemma~\ref{fg1}     for the function 
   \begin{equation}
G (\xi) =    \sqrt{2\pi}    \Gamma (1 + i \xi) (\Phi\Omega)(\xi) 
 \label{eq:LM1}\end{equation}
 plays the crucial role.

Since $\theta \in \cal Z'$, its   Fourier transform $a=\Phi \theta$ is   correctly defined as   an element  of   $C_{0}^\infty ({\Bbb R})'$. 
Formally,
    \begin{equation}
 a(\xi) =  (\Phi \theta)(\xi)=  (2\pi)^{-1/2}\int_{0}^\infty   h (t) t^{ -i\xi}   dt ,
\label{eq:M6}\end{equation}
that is, $a(\xi)$ is the Mellin transform of the function $  h (t) t^{1/2}$. Let $\Omega\in {\cal Z}$.  Passing  to the Fourier transforms and using notation \e{eq:LM1},  we   see that
\begin{equation}
  { \la} \theta, \Omega { \ra}={ \la} a,  \Phi\Omega  { \ra}={ \la} b,  G { \ra}
 \label{eq:LM3}\end{equation}
 where $G\in C_{0}^\infty ({\Bbb R})$ and the distribution $b\in C_{0}^\infty ({\Bbb R})'$ is given by the relation
\begin{equation}
 b(\xi)  =  (2\pi)^{-1/2} a(\xi)  \Gamma (1-i\xi)^{-1} 
\label{eq:M9}\end{equation}
which is of course the same as  \e{eq:Bb}. Thus we are led to the following

\begin{definition}\label{HBS}
Let $h \in  {\cal Z}'_{+}$. The distribution $b\in C_{0}^\infty ({\Bbb R})'$   defined  by formulas \e{eq:HH3}, \e{eq:M6} and \e{eq:M9} is called the $b$-function of the kernel $h(t)$ (or of the Hankel operator $H$). Its  Fourier transform $s= \sqrt{2\pi}    \Phi^* b\in {\cal Z} '$ is called the $s$-function or the sign-function.
  \end{definition}
  
  Recall that the distribution 
   $h^\natural \in  {\cal Z}_{+}'$ was defined by  relation \e{eq:B}. The  following assertion is an immediate consequence of formulas \e{eq:HH3}, \e{eq:M6} and \e{eq:M9}.
  
   \begin{proposition}\label{1wx}
 The mappings
    \[
    h\mapsto \theta\mapsto a \mapsto b \mapsto s \mapsto h^\natural 
    \]
    yield   one-to-one correspondences $($bijections$)$
    \[
    {\cal Z}_{+}' \to  {\cal Z}' \to C_{0}^\infty ({\Bbb R})  \to C_{0}^\infty ({\Bbb R})\to  {\cal Z}' \to   {\cal Z}_{+}'.
    \]
    All of them, as well as their inverse mappings, are continuous.
 \end{proposition}

     Putting together equalities  \e{eq:HH4} and \e{eq:LM3}, we see that
\begin{equation}
  {\pmb\la} h, F {\pmb\ra} = { \la} b,  G { \ra}.
 \label{eq:LL}\end{equation}
  Combining this relation with Lemma~\ref{fg1} and Definitions~\ref{HH}, \ref{HBS} and using notation \e{eq:MAID1}, we obtain the main identity \e{eq:MAID}. To be more precise, we have the following result.
 
 \begin{theorem}\label{1}
 Suppose that  $h \in  {\cal Z}'_{+}$, and let  $b\in C_{0}^\infty ({\Bbb R})'$   be the corresponding $b$-function.   
  Let $f_{j}\in {\cal D}$, $j=1,2$, and  let the functions      $ g_{j}=\Xi  f_{j}  $ be defined by formula \e{eq:M10}. Then    $ g_{j} \in C_{0}^\infty ({\Bbb R})$ and the representation 
\begin{equation}
{ \pmb\la} h ,  \bar{f}_{1}\star f_{2} {\pmb \ra} =    \la b,   ({\cal J}\bar{g}_1) * g_{2} \ra =:  b [g_{1}, g_{2}]  
\label{eq:M15s}\end{equation}
  holds.
 \end{theorem}
 
 Passing in the right-hand side of \e{eq:M15s} to the Fourier transforms and using that
 \[
 \Phi^*( ({\cal J}\bar{g}_1) * g_{2})=  (2\pi)^{1/2} \ov{\Phi^* g_1} \Phi^* g_2,    
 \]
 we obtain
 
   \begin{corollary}\label{four}
 Let $s \in {\cal Z} '$ be the sign-function of $h$, and let  $u_{j}=\Phi^* g_{j} =\Phi^* \Xi f_{j} \in {\cal Z} $. Then
 \begin{equation}
{ \pmb\la} h ,  \bar{f}_{1}\star f_{2} {\pmb \ra} =   \la s,   \bar{u}_{1} u_2\ra =:  s [u_{1}, u_{2}]  . 
\label{eq:M15f}\end{equation}
 \end{corollary}

 Loosely speaking, equalities \e{eq:M15s} and \e{eq:M15f} mean that
 \begin{equation}
 \begin{split}
{ \pmb\la} h ,  \bar{f}_{1}\star f_{2} {\pmb \ra}  =&    \int_{-\infty}^\infty\int_{-\infty}^\infty b(\xi-\eta)    g_{1}(\eta)   \overline{g_{2} (\xi)} d\xi d\eta
 \\
=& \int_{-\infty}^\infty s(x)   u_{1}(x)     \overline{u_2(x) } dx. 
  \end{split}
\label{eq:M15}\end{equation}  
    
    In the particular case $h(t)=t^{-1}$, we have $\theta (x)=1$, 
    \begin{equation}
     a(\xi)=    (2\pi )^{1/2} \d (\xi),\q b(\xi)=   \d (\xi), \q s(x)=  1,
     \label{eq:carl1}\end{equation} 
      and hence \e{eq:M15} yields
   \[
    { \pmb\la} h ,  \bar{f}_{1}\star f_{2} {\pmb \ra} =     \int_{-\infty}^\infty    g_{1}(\xi)    \overline{g_{2}(\xi)} d\xi = \int_{-\infty}^\infty |\Gamma (1/2+i \xi)|^2   \ti{f}_{1}(\xi)    \overline{ \ti{f}_{2}(\xi)} d\xi ,
    \]
   where
 \begin{equation}
    |\Gamma (1/2+i \xi)|^2 =\frac{\pi} {\cosh (\pi \xi)}.
 \label{eq:carl}\end{equation}  
    This leads to the familiar diagonalization of the Hankel operator $H$ with kernel $h(t)=t^{-1}$. This operator is known as   the Carleman operator and will be denoted by $\bf C$.

      \medskip
 
 {\bf 2.3.}
 According to Proposition~\ref{1wx} the distribution $ h^\natural$ determines uniquely the distribution $ h$.       Let us now obtain an explicit formula for the mapping $ h^\natural \mapsto h$. This requires some auxiliary information. 
   
    Let ${\pmb \Gamma}_{\alpha}: C_{0}^\infty  ({\Bbb R})\to C_{0}^\infty  ({\Bbb R})$, $\alpha>0$, be the operator of multiplication by the function $\Gamma(\alpha+i\xi)$.    Making the change of variables $t=e^{-x}$ in the definition of the gamma function, we see that
\[
\Gamma (\alpha+ i\lambda)=\int_0^\infty e^{-t} t^{\alpha + i\lambda-1}dt =\int_{-\infty}^\infty e^{-e^{-x}} e^{-\alpha x} e^{- ix\lambda}dx, \q \alpha>0,
\]
and hence  
 \begin{equation}
(2\pi)^{-1}\int_{-\infty}^\infty e^{ ix\lambda}\Gamma (\alpha+i\lambda) d\lambda=
 e^{-e^{-  x}} e^{-\alpha x} .   \label{eq:IDGg}\end{equation}
It follows that
 \begin{equation}
(\Phi^* {\pmb \Gamma}_{\alpha} \Phi \Omega)(x)= \int_{-\infty}^\infty e^{\alpha(y- x)} e^{-e^{y-  x}} \Omega(y)dy.
\label{eq:IDGgk}\end{equation}

   Let us also introduce the   operator ${\sf L}_{\alpha}$:
   \[
( {\sf L}_{\alpha}F)(\lambda)= \lambda^\alpha \int_{0}^\infty e^{-t\lambda} t^{\alpha-1} F(t) dt, \q \lambda>0, \q \alpha>0.
\]
 Obviously, ${\sf L}_{\alpha}F \in C^\infty ({\Bbb R}_{+})$
for  all bounded functions $F(t)$ and, in particular, for $F \in {\cal Z}_{+}$.
Note that ${\sf L}_{\alpha}$ is the Laplace operator  $\sf L$, 
\begin{equation}
  ({\sf L} F) (\lambda)=\int_{0}^\infty e^{-t \lambda } F(t) dt,  
\label{eq:LAPj}\end{equation}
sandwiched  by the weights $ \lambda^\alpha$ and $ t^{\alpha-1}$. 

Recall that the operator $\cal R$ defined by \e{eq:HH3s}  is a one-to-one mapping of $ {\cal Z}_{+}$ onto $ \cal Z $.

  We need the following result.

  \begin{lemma}\label{LAPL}
 For all $\alpha>0$,  the identity
    \begin{equation}
{\sf L}_{\alpha}= {\cal R}^{-1}{\cal J} \Phi^* {\pmb \Gamma}_{\alpha} \Phi {\cal R}
 \label{eq:LAPL1}\end{equation}
 holds. In  particular, ${\sf L}_{\alpha} $ as well as its inverse are the one-to-one continuous mappings  of $ {\cal Z}_{+}$ onto itself.
  \end{lemma}
  
   \begin{pf}
  Putting $\Omega(y)= ({\cal R} F)(y)= F(e^y)$ in \e{eq:IDGgk} and making the change of variables $t=e^y$, we find that
\[
(\Phi^* {\pmb \Gamma}_{\alpha} \Phi {\cal R} F)(x)= e^{-\alpha  x}\int_{0}^\infty   e^{-e^{-  x}t}t^{\alpha-1} F(t)dt.
\]
Now making the change of variables $\lambda=e^{-x}$,  we arrive at the identity \e{eq:LAPL1}.

Consider the right-hand side of  \e{eq:LAPL1}. All mappings  $ {\cal R} :  {\cal Z}_{+}\to  {\cal Z}$,
$\Phi:  {\cal Z}\to  C_{0}^\infty ({\Bbb R} )$,  
$ {\pmb\Gamma}_{\alpha}  : C_{0}^\infty ({\Bbb R} )\to  C_{0}^\infty ({\Bbb R} )$,  $\Phi^*:     C_{0}^\infty ({\Bbb R} )\to {\cal Z} $,   ${\cal J}:  {\cal Z}\to  {\cal Z}$ are bijections. All of them as well as their inverses are continuous.  Therefore  the identity \e{eq:LAPL1} ensures the same result for the operator   ${\sf L}_{\alpha}: {\cal Z}_{+}\to {\cal Z}_{+}$.
    \end{pf}

  To recover $h(t)$, we proceed from formula \e{eq:LL}.
    Passing to the Fourier transforms, we can write it as
     \[
  {\pmb\la} h, F {\pmb\ra} = (2\pi)^{-1/2} { \la} s,  \Phi^*G { \ra}
\]
 where $G$ is defined by formulas \e{eq:HH3s}, \e{eq:LM1}, that is,
 $
 G= (2\pi)^{1/2}{\pmb\Gamma}_{1} \Phi {\cal R} F.
$ 
Therefore using the identity  \e{eq:LAPL1}  for $\alpha=1$, we see that
 \[
  {\pmb\la} h, F {\pmb\ra} = { \la} s,  {\cal J} {\cal R} {\sf L}_{1} F{ \ra}= \int_{-\infty}^\infty s(x)\overline{  ({\sf L}_{1} F)(e^{- x})}  dx.
\]
   Making   the change of variables $\lambda =e^{-x}$ in the right-hand side, 
we obtain  the identity
\[
 {\pmb\la} h, F {\pmb\ra}=   {\pmb\la}  h^\natural , {\sf L}_{1} F {\pmb\ra}.
\] 
  Passing here  to adjoint operators  and using that $F\in  {\cal Z}_{+} $ is arbitrary, we find that
 \begin{equation}
h= {\sf L}_{1}^* h^\natural ,
\label{eq:T2x}\end{equation}
which gives the precise sense to formula \e{eq:conv1}.  Of course formula  \e{eq:conv1} can also be rewritten as
 \[    
  h (t)  =    \int_{-\infty}^\infty     e^{-t e^{-x}} e^{-x} s (x)  d x.
  \]

    Let us state the result obtained. 
    
    \begin{theorem}\label{round}
  Let    $h \in  {\cal Z}_{+}'$, and let     $ s\in {\cal Z}'$  be the corresponding sign-functions 
  $($see Definition~\ref{HBS}$)$. Define
 the distribution $ h^\natural$   by formula   \e{eq:B}. Then $  h^\natural (\lambda)$ belongs to the class $  {\cal Z}_{+}'$ and $h$ can be recovered from $ h^\natural$  by formula \e{eq:T2x}.
  \end{theorem}
        
    We emphasize  that in  the roundabout
  $  h\mapsto  h^\natural \mapsto h$
    the mappings $h\mapsto h^\natural $ as well as its inverse  $ h^\natural \mapsto h$ are one-to-one continuous mappings of the set $ {\cal Z}_{+}'$ onto itself.

  Let us also give   a direct expression of $u (x)=(\Phi^* g)(x)$ in terms of $f(t)$.

  \begin{lemma}\label{LM}
  Suppose that $f\in  {\cal D}$ and put $ \varphi(t)=t^{1/2} f(t)$.
  Let    $g(\xi)$   be defined by formula \e{eq:MAID1} and   $u (x) =(\Phi^* g)(x) $. Then 
  \begin{equation}
u (x)  =   ({\sf L}_{1/2}  \varphi) (e^{-x}).
\label{eq:LAPg}\end{equation}
  \end{lemma}
  
   \begin{pf}
  Since $({\cal R} \varphi)(x)= (Uf) (x)$, 
   it follows from  formula \e{eq:LAPL1} for $\alpha=1/2$ that
      \[
  ( {\cal R}^{-1}{\cal J} \Phi^* {\pmb \Gamma}_{1/2} \Phi U f)(\lambda)= ({\sf L}_{1/2} \varphi)(\lambda)  .
 \]
The left-hand side here equals $  ( {\cal R}^{-1}{\cal J}u)(\lambda)$ which after the change of variables $\lambda=e^{-x}$ yields \e{eq:LAPg}.
   \end{pf}
   
Now   we can rewrite identity \e{eq:M15f} in a slightly different way.
   
     \begin{corollary}\label{LMC}
  Let   $h\in  {\cal Z}_{+}'$, and let the distribution $h^\natural\in  {\cal Z}_{+}'$ be defined by formula \e{eq:B}. Then for arbitrary $ f_j\in  {\cal D}$, $j=1,2$, and $\varphi_{j}(t)=t^{1/2} f_{j}(t)$, we have
     \begin{equation}
{ \pmb\la} h ,  \bar{f}_{1}\star f_{2} {\pmb \ra} =  { \pmb\la} h^\natural , \overline{ {\sf L}_{1/2}  \varphi_1} {\sf L}_{1/2}  \varphi_{2} {\pmb \ra} .
\label{eq:M15g}\end{equation} 
  \end{corollary}
  
   \begin{pf}
   It suffices to make the change of variables $x=-\ln \lambda$ in the right-hand side of  \e{eq:M15} and to take equality \e{eq:LAPg} into account.  
   \end{pf}
  
  We emphasize that according to Lemma~\ref{LAPL},  ${\sf L}_{1/2} \varphi_{j} \in  {\cal Z}_{+}$  and hence  $ \overline{ {\sf L}_{1/2}  \varphi_1} {\sf L}_{1/2}  \varphi_{2}  \in  {\cal Z}_{+}$. Thus  the right-hand side of \e{eq:M15g} is correctly defined.

     \medskip
 
 {\bf 2.4.}
 Finally, we check that a distribution $h\in   {\cal Z}_{+}'$ is determined uniquely by the values $ { \pmb\la} h ,  \bar{f}_{1}\star f_{2} {\pmb \ra}$ on $f_1, f_2 \in \cal D$. First we consider convolution operators. Let us introduce the shift in the space $L^2 ({\Bbb R})$:
    \begin{equation}
(T (\tau)g)(\xi)= g(\xi-\tau), \q  \tau\in {\Bbb R}.
  \label{eq:shift}\end{equation}
  Since
  \[
(  g_{1}* g_{2}) (\xi)=\int_{-\infty}^\infty (T(\tau) g_{1} ) (\xi) g_{2}(\tau)d\tau,\q \forall g_{1},  g_{2}\in C_{0}^\infty ({\Bbb R}),
  \]
  we have the formula
    \begin{equation}
   \la b,  ( {\cal J} \bar{g}_{1})* g_2\ra =   \int_{-\infty}^\infty  \la b,  T (\tau) {\cal J} \bar{g}_{1}\ra   \ov{{g}_{2}(\tau)} d \tau
 \label{eq:CF} \end{equation}
 where for $b\in C_{0}^\infty ({\Bbb R})'$ the function $ \la b,  T (\tau) {\cal J} \bar{g}_{1}\ra $ is  infinitely differentiable in $\tau\in {\Bbb R}$. 
 
   The following assertion is quite standard.

\begin{lemma}\label{CC}
 Let $b\in C_{0}^\infty ({\Bbb R})'$. Suppose that   $\la b,  ( {\cal J} \bar{g}_{1})* g_2\ra = 0$ for all      $ g_1, g_{2} \in C_{0}^\infty ({\Bbb R})$.  Then $b=0$.
   \end{lemma}
   
   \begin{pf}
 If $\la b,  ( {\cal J} \bar{g}_{1})* g_2\ra = 0$ for all      $ g_2  \in C_{0}^\infty ({\Bbb R})$, then $\la b,  T (\tau) {\cal J} \bar{g}_{1}\ra  =0$ for all $\tau\in{\Bbb R}$ according to formula \e{eq:CF}. In particular, for $\tau=0$ we have $\la b,    {\cal J} \bar{g}_{1}\ra  =0$ whence  $b=0$ because  $g_1 \in C_{0}^\infty ({\Bbb R})$
   is arbitrary.
   \end{pf}
   
   Next we pass to Hankel operators.
   
   \begin{proposition}\label{CC1}
 Let $h\in  {\cal Z}'_{+}$. Suppose that   ${ \pmb\la} h ,  \bar{f}_{1}\star f_{2} {\pmb \ra} = 0$ for all      $ f_1, f_{2} \in  \cal D$.  Then $h=0$.
   \end{proposition}
   
   \begin{pf}
  Let $b\in C_{0}^\infty ({\Bbb R})'$ be the $b$-function of $h$ (see Definition~\ref{HBS}). For arbitrary 
    $ g_1 , g_{2} \in C_{0}^\infty ({\Bbb R})$, we can construct $ f_1, f_{2} \in  \cal D$ by formula \e{eq:M10}.  Since ${ \pmb\la} h ,  \bar{f}_{1}\star f_{2} {\pmb \ra} = 0$, it follows from  the identity \e{eq:M15s} that  $\la b,  ( {\cal J} \bar{g}_{1})* g_2\ra = 0$. Therefore  $b=0$ according to Lemma~\ref{CC}. Now Proposition~\ref{1wx} implies that $h=0$.
      \end{pf} 
      
        %%%%%%%%%%%%%%%%%%%%%%%%%%%
  
\section{Bounded Hankel operators}  

%%%%%%%%%%%%%%%%%%%%%%%%%%

Our main goal here is to show that the condition $h\in  {\cal Z}'_{+}$ is satisfied for all bounded Hankel operators $H$.

\medskip
  
  {\bf 3.1.}
   In this section we a priori only assume that $h\in C_{0}^\infty ({\Bbb R_{+}})'$ and consider the Hankel form \e{eq:HH2} on functions $f_{1}, f_{2} \in C_{0}^\infty ({\Bbb R_{+}})$. Let $T_{+}(\tau)$ where $\tau \geq 0$ be the restriction of the shift \e{eq:shift} on its invariant subspace $L^2 ({\Bbb R}_{+})$. Since
   \[
   (\bar{f}_{1} \star f_2)(t) =\int_{0}^\infty (T_{+}(\tau)\bar{f}_{1}) (t) f_2(\tau)d\tau , \q \forall f_1, f_{2} \in  C_{0}^\infty ({\Bbb R_{+}}),
   \]
   for all $h\in C_{0}^\infty ({\Bbb R_{+}})'$ we have the formula
    \begin{equation}
  {\pmb \la } h,  \bar{f}_{1} \star f_2  {\pmb \ra } =   \int_{0}^\infty   {\pmb \la } h, T_{+}(\tau) \bar{f}_{1}   {\pmb \ra }   \ov{f_{2}(\tau)} d \tau .
 \label{eq:CFH} \end{equation}
Here   the function $  {\pmb \la } h, T_{+}(\tau) \bar{f}_{1}   {\pmb \ra } $ is  infinitely differentiable in $\tau\in {\Bbb R}_{+}$, and this function, as well as all its derivatives, have finite limits as $\tau\to 0$. 
    In the theory of Hankel operators,  formula \e{eq:CFH} plays  the role of formula \e{eq:CF} for convolution operators.   
    
  The proof of the  following assertion  is almost the same as  that of Lemma~\ref{CC}. 
   
   \begin{proposition}\label{CC2}
 Let $h\in C_{0}^\infty ({\Bbb R_{+}})'$. Suppose that   ${ \pmb\la} h ,  \bar{f}_{1}\star f_{2} {\pmb \ra} = 0$ for all      $ f_1, f_{2} \in  C_{0}^\infty ({\Bbb R_{+}})$.  Then $h=0$.
   \end{proposition}
   
      \begin{pf}
       If $  {\pmb \la } h,  \bar{f}_{1} \star f_2  {\pmb \ra } = 0$ for all      $ f_2  \in C_{0}^\infty ({\Bbb R}_{+})$, then $ {\pmb \la } h, T_{+}(\tau) \bar{f}_{1}   {\pmb \ra }  =0$ for all $\tau\in [0,\infty)$ according to formula \e{eq:CFH}. In particular, for $\tau=0$ we have $ {\pmb \la } h,  \bar{f}_{1}   {\pmb \ra }  =0$  which implies that $h=0$ because  $f_1 \in C_{0}^\infty ({\Bbb R}_{+})$   is arbitrary.
   \end{pf}
   
  Of course Propositions~\ref{CC1} and \ref{CC2}  differ only by   the set of functions on which the Hankel form is considered.

   Assume now that 
    \begin{equation}
 |  {\pmb \la } h, \bar{f}\star f  {\pmb \ra } |\leq C  \| f\|^2,\q \forall f\in C_{0}^\infty ({\Bbb R_{+}}).
  \label{eq:bbb1}\end{equation} 
  Then there exists a bounded   operator $H$ such that 
  \begin{equation}
  (H f_{1},f_{2})=  {\pmb \la } h, \bar{f}_{1}\star f_{2} {\pmb \ra } ,\q \forall f_{1},f_{2} \in C_{0}^\infty ({\Bbb R_{+}}).
  \label{eq:Bbb}\end{equation} 
{\it  We call $H$ the Hankel operator associated to the Hankel form }  ${\pmb \la } h, \bar{f}_{1}\star f_{2} {\pmb \ra }$. 
  
 \medskip
  
  {\bf 3.2.}
  It is possible to characterize Hankel operators by some commutation relations. Let us define a bounded operator $\Sigma$ in the space $ L^2({\Bbb R}_{+})  $ by the equality
    \[
  (\Sigma f) (t)= -2 e^{-t}\int_{0}^t e^s f(s)ds.
 \]
  Note that
    \begin{equation}
 \Sigma= -2 \int_{0}^\infty T_{+}(\tau ) e^{-\tau} d\tau.
  \label{eq:CR2x}\end{equation}
  
    \begin{lemma}\label{CC3}
 Let assumption \e{eq:bbb1} hold. Then the corresponding Hankel operator $H$ satisfies the commutation relations:
    \begin{equation}
H T_{+}(\tau) =T_{+}(\tau) ^* H,\q \forall \tau\geq 0,
  \label{eq:CR1}\end{equation}
  and  
     \begin{equation}
 H\Sigma =\Sigma^* H.
  \label{eq:CR2}\end{equation}
   \end{lemma}
   
 \begin{pf}
Since
 \[
 (T_{+}(\tau)\bar{f}_{1})\star f_{2} =  \bar{f}_{1}\star (T_{+}(\tau) f_{2}),\q \forall \tau\geq 0,
 \]
relation  \e{eq:CR1} directly follows from definition  \e{eq:Bbb}. By virtue of formula \e{eq:CR2x},
 relation \e{eq:CR2} is a consequence of \e{eq:CR1}.
   \end{pf}
   
   Below we need the   Nehari theorem; see the original paper \cite{Nehari}, or the books 
   \cite{Pe}, Chapter~1, \S 1 or \cite{Po}, Chapter~1, \S 2. We formulate it     in      the Hardy space  ${\tt H}_{+}^2({\Bbb R})$  of functions analytic in the upper half-plane. We denote by $\wh{\Sigma}$   the operator of multiplication by the function $(\mu-i)/(\mu+i)$ in      this  space.

          \begin{theorem}[Nehari]\label{Neh}
          Let $\omega\in L ^\infty({\Bbb R})$, and let a bounded operator $\wh{H}$ in  the space    ${\tt H}_{+}^2({\Bbb R})$ be defined by the relation
             \begin{equation}
 (\wh{H} \hat{f}_{1}, \hat{f}_2)=\int_{- \infty}^\infty \omega(\mu) \hat{f}_{1}(-\mu) \ov{ \hat{f}_2(\mu)} d \mu, \q \forall \hat{f}_{1}, \hat{f}_{2}\in {\tt H}_{+}^2({\Bbb R}).
  \label{eq:Neh}\end{equation}
  Then $\wh{H}\wh{\Sigma}=\wh{\Sigma}^*\wh{H}$. Conversely, 
       if $\wh{H}$ is a bounded operator in      ${\tt H}_{+}^2({\Bbb R})$ and $\wh{H}\wh{\Sigma}=\wh{\Sigma}^*\wh{H}$, then there exists a function $\omega\in L ^\infty({\Bbb R})$ such that representation \e{eq:Neh} holds.
         \end{theorem}
   
   The following assertion can be regarded as a   translation of this theorem into the space $L^2 ({\Bbb R}_{+})$. Recall that the Fourier transform $\Phi : {\tt H}_{+}^2({\Bbb R})\to L^2 ({\Bbb R}_{+})$ is the unitary operator. Since
 \[
 \int_{-\infty}^\infty \frac{\mu-i}{\mu+i}e^{-i \mu t} d \mu= -4\pi e^{-t}
 \]
 for $t>0$ and this integral is zero for $t<0$,   we have the relation
    \begin{equation}
    \wh{\Sigma}=\Phi^* \Sigma \Phi.
     \label{eq:sigma}\end{equation}
   
       \begin{theorem}\label{NehC}
       If representation
    \e{eq:Bbb} holds with $h=\Phi \omega$ where $\omega\in L^\infty ({\Bbb R})$ $($in this case $h  \in {\cal S}'\subset C_{0}^\infty ({\Bbb R}_{+})')$, then estimate \e{eq:bbb1} is true and the corresponding Hankel operator satisfies commutation relation  \e{eq:CR1}. Conversely, 
 if a bounded operator $H$ satisfies    \e{eq:CR1}, then representation
    \e{eq:Bbb} holds with $h=\Phi \omega$ for some $\omega\in L^\infty ({\Bbb R})$. 
   \end{theorem}
   
    \begin{pf}   
    Since
      \[
(\Phi^* (f_{1} \star \bar{f}_{2}))(\mu)= ( {\cal J} \hat{f}_{1})(\mu) \ov{ \hat{f}_2(\mu)}, \q \forall f_1, f_{2}\in C_{0}^\infty({\Bbb R}_{+}),  
\]
 where $   \hat{f}_{1}=\Phi^* f_{1},  \hat{f}_2=\Phi^* f_2$ and $ ( {\cal J} \hat{f}_{1})(\mu) = \hat{f}_{1}(-\mu)$,    we have
    \begin{equation}
 {\pmb\la }h, \bar{f}_{1}\star f_{2} {\pmb\ra } = {\pmb\la }\Phi^* h, \ov{( {\cal J} \hat{f}_{1})}
\, \hat{f}_{2} {\pmb\ra }, \q \forall h\in {\cal S}'.
      \label{eq:Neh2}\end{equation}
      Therefore estimate   \e{eq:bbb1} is satisfied  if $\Phi^*h\in L^\infty ({\Bbb R})$. Relation  \e{eq:CR1} for the corresponding Hankel operator  $H$ follows from Lemma~\ref{CC3}.
      
      Conversely, 
 if a bounded operator $H$ satisfies relation  \e{eq:CR1}, then by virtue of  \e{eq:CR2x} it also satisfies relation  \e{eq:CR2}. Hence it follows from \e{eq:sigma} that $\wh{H}\wh{\Sigma}=\wh{\Sigma}^*\wh{H}$ where  $\wh{H}=\Phi^* H \Phi$ is a bounded operator  in      the  space  ${\tt H}_{+}^2({\Bbb R})$.  Thus, by Theorem~\ref{Neh},   there exists a function $\omega\in L ^\infty({\Bbb R})$ such that representation \e{eq:Neh} holds. It means that
    \begin{equation}
 (H f_{1}, f_2)=\int_{- \infty}^\infty \omega(\mu) \hat{f}_{1}(- \mu) \ov{ \hat{f}_2(\mu)} d\mu, \q \forall f_{1}, f_{2}\in L^2({\Bbb R}_{+}).
  \label{eq:Neh3}\end{equation}
  If $h=\Phi \omega$, then the right-hand sides in \e{eq:Neh2} and \e{eq:Neh3} coincide.
This yields representation    \e{eq:Bbb}.
       \end{pf}
       
        \begin{corollary}\label{NehC1}
 For a bounded operator $H$,     commutation relations  \e{eq:CR1} and  \e{eq:CR2} are equivalent. 
   \end{corollary}
   
    \begin{pf}
    As was already noted,    \e{eq:CR2} follows from \e{eq:CR1} according to formula \e{eq:CR2x}. Conversely, if $H$ satisfies \e{eq:CR1}, then representation    \e{eq:Bbb} holds according to Theorem~\ref{NehC}. Thus it remains to use Lemma~\ref{CC3}.
       \end{pf}

  {\it A function $\omega \in L^\infty   ({\Bbb R})$ such that $\Phi\omega  =h $    is called the   symbol of 	 a bounded Hankel operator  $H$ with kernel} $h(t)$. Of course if $\omega\in {\sf H}^\infty_{-}   ({\Bbb R})$, that is, $\omega$ admits an analytic continuation to a bounded function in the lower half-plane,  then the corresponding Hankel operator is zero. Therefore the symbol is defined up to a function in the class ${\tt H}^\infty_{-}   ({\Bbb R})$.  
   
       \medskip
  
  {\bf 3.3.}
  Now we are in a position to check that the condition $h\in  {\cal Z}'_{+}$ is satisfied for all bounded Hankel operators. By definition  \e{eq:heta} it means that   distribution  \e{eq:HH3}  belongs to the class $ {\cal Z}'$. We shall verify the stronger inclusion $\theta\in {\cal S}'$.
  
  To that end, it suffices to check that, for some $N\in {\Bbb Z}_{+}$ and some $\kappa\in {\Bbb R}$,
  \begin{equation}
| \la \theta,\Omega\ra |\leq C \sum_{n=0}^N \max_{x\in{\Bbb R}} \big( (1+|  x|)^\kappa  |\Omega^{(n)} (x)|\big), \q \forall \Omega \in C_{0}^\infty ({\Bbb R} ).
  \label{eq:ch}\end{equation}
 Putting $F(t)=\Omega (\ln t)$, we see that \e{eq:ch} is equivalent to the estimate
   \begin{equation}
| \la h, F\ra| \leq C\sum_{n=0}^N \max_{t\in{\Bbb R}_{+}} \big( (1+| \ln t|)^\kappa  t^n |F^{(n)} (t)|\big), \q F\in C_{0}^\infty ({\Bbb R}_{+}). 
  \label{eq:CHb}\end{equation}
  
  If $h\in L^1_{\rm loc}({\Bbb R})$, then 
    estimate \e{eq:CHb}  for $N=0$ is equivalent to the convergence of integral  \e{eq:ass}. If $H$ is Hilbert-Schmidt, that is
  \[
  \int_{0}^\infty |h(t)|^2 tdt<\infty, 
  \]
 then integral  \e{eq:ass} converges for any $\kappa >1/2$.  Similarly, if $|h(t)| \leq C t^{-1}$, then integral  \e{eq:ass} converges for any $\kappa >1$.    
  
 For the proof of \e{eq:CHb} in the general case, we use the following elementary result. Its proof is given in Appendix~A.
  
    \begin{lemma}\label{CHa}
If  $F\in C_{0}^\infty ({\Bbb R}_{+})$, then for an arbitrary    $\kappa>2$ the estimate
 \[
\| \Phi^* F\| _{L^1({\Bbb R})}\leq C\sum_{n=0}^2 \max_{t\in{\Bbb R}_{+}} \big((1+| \ln t|)^\kappa t^n |F^{(n)} (t)|\big)
 \]
  holds.
   \end{lemma}

   \begin{corollary}\label{CH}
If $h= \Phi\omega$ where $\omega\in L^\infty ({\Bbb R})$, then estimate \e{eq:CHb} holds for $N=2$ and  an arbitrary $\kappa>2$.  
   \end{corollary}
   
   Since, by Theorem~\ref{NehC}, for a bounded Hankel operator $H$, its kernel $h= \Phi\omega$ for some $\omega\in L^\infty ({\Bbb R})$, we arrive at the following result.
   
    \begin{theorem}\label{CHb}
Suppose that $h\in C_{0}^\infty ({\Bbb R}_{+})'$ and that condition \e{eq:bbb1} is satisfied.  Then estimate \e{eq:CHb} holds for $N=2$ and  an arbitrary $\kappa>2$; in particular, $h\in  {\cal Z}'_{+}$.
   \end{theorem}

    The following simple example shows that   for $N=0$ estimate \e{eq:CHb} is not in general true (for all $\kappa$).  
   
   \begin{example}\label{ECH1}
   Let $h(t)=e^{-i t^2}$. Then the corresponding Hankel operator $H$ is bounded because
 according to the formula
   $
  e^{- i (t+s)^2}=  e^{-i t^2}e^{-i 2ts}e^{-i s^2},
$
  it is a product of three bounded operators. However integral  \e{eq:ass} diverges  at infinity  for all $\kappa$.    In this example   condition \e{eq:CHb} is satisfied for $N=1$ and   $\kappa=0$. Indeed, integrating by parts, we see that
  \begin{equation}
  \int_{0}^\infty h(t) \ov{F(t)} dt=-   \int_0^\infty h_{1}(t) \ov{F'(t)} dt
   \label{eq:CHbk}\end{equation}
  where the function
 $  h_{1}(t)=\int_{0}^t e^{-i s^2} ds$ is bounded and   $h_{1}(t)= O( t )$
  as $t\to 0$. Therefore integral \e{eq:CHbk} is bounded by $\max_{t\in{\Bbb R}_{+}} \big(  (1+|\ln t|)^\kappa t |F' (t)|\big)$ for $\kappa>1$.
   \end{example}

  Note that  the symbol of $H$ equals $\omega (\mu)= ( 2\pi)^{-1} e^{-\pi i/4} e^{i\mu^2/4}$.  More generally, one can consider the class of symbols $\omega (\mu)$ such that $\omega\in C^\infty ({\Bbb R})$, $\omega (\mu)= e^{i \omega_{0} \mu^\alpha}$,    $\omega_{0}> 0$, for large positive $\mu$ and $\omega (\mu)= 0$ for large negative $\mu$. Of course  Hankel operators    with such symbols are bounded.  Using the stationary phase method, we find that for $\alpha>1$ the corresponding kernel $h(t)$   has the asymptotics
   \begin{equation}
  h(t)   \sim h_{0} t^{\beta} e ^{i \sigma t^{\gamma}}, \q t\to \infty,
     \label{eq:CHbk1}\end{equation}
  where $\beta= (1-\alpha/2)(\alpha-1)^{-1}$, $\gamma= \alpha (\alpha-1)^{-1}$ and $h_{0}$, $\sigma=\bar{\sigma}$ are some  numbers. Moreover, $h(t)$ is a bounded function on all finite intervals. Similarly to Example~\ref{ECH1}, it can be checked that for such kernels  condition \e{eq:CHb} is satisfied for $N=1$ but not for  $N=0$. If $\alpha\in (0,1)$, then $h(t)$   has   asymptotics \e{eq:CHbk1}    for $t\to 0$.

       \medskip
  
  {\bf 3.4.}
  Here we shall show that, for bounded Hankel operators $H$,  the representations \e{eq:M15s} and \e{eq:M15f} extend to all $f_{1}, f_{2}\in L^2({\Bbb R}_{+})$. By Theorem~\ref{CHb}, we have    $h \in  {\cal Z}_{+}'$. Let $b$ and $s$ be the corresponding $b$- and $s$-functions (see Definition~\ref{HBS}). 
 Recall that the operator $ \Xi: L^2({\Bbb R}_{+})\to L^2({\Bbb R} )$ was defined by formula \e{eq:MAID1}. We denote by $K$   the operator of multiplication by the function $\sqrt{\cosh(\pi\xi)/\pi}$ in the space $L^2({\Bbb R} )$. It follows from identity \e{eq:carl}  and the unitarity of the Mellin transform \e{eq:M1}  that 
 \[
 \| K \Xi f\|= \| f\|,
\]
 and hence the operator $K \Xi: L^2({\Bbb R}_{+})\to L^2({\Bbb R} )$  is unitary. Therefore in view of the identities  \e{eq:M15s} and  \e{eq:M15f} we have the following result.
  
    \begin{lemma}\label{bbb}
The  inequalities     
  \e{eq:bbb1}, 
%  for all $f\in {\cal D}$,
   \begin{equation}
 | \la b, ({\cal J}\bar{g})* g\ra |\leq C  \| Kg \|^2,\q \forall g\in C_{0}^\infty ({\Bbb R} ),
  \label{eq:bbb2}\end{equation}
  and
   \begin{equation}
 | \la s, |u|^2 \ra |\leq C  \| K\Phi u\|^2,\q \forall u\in {\cal Z},
  \label{eq:bbb3}\end{equation}
  are equivalent.  The   Hankel operator corresponding  to form  \e{eq:HH2}  is bounded if and only if one of equivalent estimates \e{eq:bbb1}, \e{eq:bbb2} or  \e{eq:bbb3} is satisfied.
 \end{lemma}
 
  These estimates  can be formulated in a slightly different way.
   Let us introduce the space ${\cal E}\subset  L^2({\Bbb R} )$ of exponentially decaying functions with the norm $\| g\|_{\cal E}=\|K g\|$. Then the space 
   ${\cal W}=\Phi^* {\cal E}$ consists of functions $u(x)$ admitting the analytic continuation $u(z)$  in the strip $\Im z\in (-\pi/2, \pi/2)$; moreover, functions  $u(x+iy)$ have limits in $ L^2 ({\Bbb R})$ as $y\to \pm \pi/2$. The identity 
   \[
   \| \Phi u\|^2_{\cal E} = (2\pi)^{-1}\int_{-\infty}^\infty \big( |u(x+i\pi/2)|^2+ |u(x-i\pi/2)|^2\big) dx= :  \| u\|^2_{\cal W}  
   \]
   defines the Hilbert norm on ${\cal W}$. We call ${\cal W}$ the exponential Sobolev space because  it is contained in standard Sobolev spaces $ {\sf H}^l ({\Bbb R}) $ for all $l$. The operators $  \Xi: L^2({\Bbb R}_{+})\to {\cal E}$ and $ \wh{\Xi}:=\Phi^* \Xi: L^2({\Bbb R}_{+})\to {\cal W}$ are of course unitary. Obviously, $\| Kg \|$ and $\| K\Phi u \|$ in the right-hand sides of \e{eq:bbb2} and \e{eq:bbb3} can be replaced by $\| g\|_{\cal E} $ and $\| u\|_{\cal W} $, respectively. Note that the inclusions $f\in L^2({\Bbb R}_{+})$,  $g=\Xi f\in {\cal E}$ and $u=\wh{\Xi} f \in {\cal W}$  are equivalent.   
   
   If one of the equivalent estimates  \e{eq:bbb1},  \e{eq:bbb2} or  \e{eq:bbb3} is satisfied, then all operators $H: L^2({\Bbb R} )\to L^2({\Bbb R} )$, $B: {\cal E} \to {\cal E}'$ and $S: {\cal W} \to {\cal W}'$ are bounded. Using that relations $f_{n}\to f$ in $L^2({\Bbb R}_{+})$, $g_{n}=\Xi f_{n}\to g =\Xi f$ in $ {\cal E}$ and $u_{n}=\Phi^* g_{n}\to u =\Phi^* g$ in $ {\cal W}$ are equivalent, we extend \e{eq:M15s} and \e{eq:M15f} to all $f\in L^2({\Bbb R}_{+})$. Thus we have obtained the following result.
 
  \begin{proposition}\label{ext}
 If one of equivalent estimates \e{eq:bbb1}, \e{eq:bbb2} or  \e{eq:bbb3} is satisfied, then
 the  identities 
   \[
 (H f_{1},f_{2})=(B g_{1},g_{2}) =(Su_{1},u_{2}), \q g_{j} =\Xi f_{j}, \q u_{j} = \Phi^* g_{j} ,
\]
%\label{eq:MAIDn}\end{equation} 
  are true for  all $f_{1}, f_{2}\in L^2({\Bbb R}_{+})$.  
   \end{proposition}

   % In particular, $H$ is compact if and only if $B_{0}$ is compact.

  Let   $ K_{l}$ be the operator of multiplication by  the function $(1+\xi^2)^{l/2}$. Then estimates \e{eq:bbb2} or  \e{eq:bbb3} are satisfied provided 
     \begin{equation}
    |  \la b, ({\cal J}\bar{g})* g\ra |\leq C_{l}  \| K_{l} g \|^2
\q {\rm or} \q 
| \la s, |u|^2 \ra |\leq C_{l}  \|   u\|^2_{{\sf H}^{l}({\Bbb R})} , 
 \label{eq:bbb5}\end{equation}
 for some $l$; in this case
  \[
  C =C_{l} \pi \max_{\xi\in{\Bbb R}}\big( (1+\xi^2)^l (\cosh(\pi\xi))^{-1}\big).
  \]
  
  We note the following assertion.
  
    \begin{proposition}\label{extb}
A Hankel operator  $H$ is bounded if its sign-function satisfies the condition
      \begin{equation}
  s\in L^1 (\Bbb R)+ L^\infty (\Bbb R).
\label{eq:LLs}\end{equation}  
 If $s \in L^\infty ({\Bbb R})$ and $s ( x)\to 0$ as $|x|\to \infty$, then  $H$   is  compact. 
   \end{proposition}
   
    \begin{pf}
    The first statement is obvious because under assumption \e{eq:LLs} the second estimate \e{eq:bbb5} is satisfied for $l>1/2$. To prove the second statement, we observe
    that  the operator $S\Phi^* K^{-1}$ is compact because  both $S$ and $K^{-1}$ are operators of multiplication by   bounded functions which tends to zero at infinity. It follows that the operator  
      $$
  H =\Xi^* \Phi (S\Phi^* K^{-1}) (K \Xi )
  $$ 
    is  also compact.  
        \end{pf}
  
 We emphasize, however, that as show already examples of Hankel operators $H$ of finite rank (see formula \e{eq:bbrr} for $h(t)=e^{-\alpha t}$), condition \e{eq:LLs} is {\it not} necessary for the boundedness of $H$.

     %************************************************************
\section{ Criteria of the sign-definiteness  }  
%*

In this section we suppose that $h(t)=\overline{h(t)}$ so that
   the  operator $H$ is symmetric.    The results of   Section~2 allow us to give simple necessary and sufficient conditions for a Hankel operator $H$ to be positive or negative. Moreover, they provide also convenient tools for a evaluation of the total multiplicity of the negative and positive spectra of $H$. We often formulate results only for the negative   spectrum.  The corresponding results   for the positive   spectrum are obtained if $H$ is replaced   by $-H$.

 \medskip
 
  {\bf 4.1.}
  Actually, we consider the problem in terms of Hankel quadratic forms rather than Hankel operators. This is both more general and more convenient.  As usual, we suppose  that  a distribution  $h \in  {\cal Z}_{+}'$ and introduce the b-function
       $b\in C_{0}^\infty ({\Bbb R})'$ and the  $s$-function  $s  \in {\cal Z}'$   by Definition~\ref{HBS}.

 Below we use the following natural notation.  
     Let ${\sf h}[\varphi , \varphi ]$ be   a real quadratic form defined on  a linear set ${\sf D} $.  We denote by $N_{\pm}({\sf h})$ the maximal dimension of linear sets ${\cal M}_{\pm}\subset {\sf D}$   such that $\pm h[\varphi,\varphi] > 0$   for all $\varphi\in {\cal  M}_{\pm}$, $\varphi\neq 0$.
  We apply this definition to the forms  $h[f ,f ]={\pmb\la}h,\bar{f} \star f {\pmb\ra}$ defined on ${\cal D} $, to $b[g ,g ]={ \la} b, {\cal J}\bar{g} * g{\ra}$ defined on $C_{0}^\infty ({\Bbb R})$ and to   $s[u,u ]={ \la} s , |u|^2{\ra}$ defined on ${\cal Z}$.  Of course, if ${\sf D} $ is dense in a Hilbert space $\cal H$ and ${\sf h}[\varphi , \varphi ]$ is closed on ${\sf D} $, then for the self-adjoint operator ${\sf  H}$ corresponding to ${\sf h}$, we have  $ N_{\pm}({\sf  H})=N_{\pm}({\sf h}) $ .
    
    Observe that formula \e{eq:MAID1} establishes one-to-one correspondence between the sets $  \cal D$ and $  C_{0}^\infty({\Bbb R})$. Moreover,   the  Fourier transform establishes one-to-one correspondence between the sets $  C_{0}^\infty({\Bbb R})$ and $  {\cal Z}  $. 
  Therefore the following assertion is a direct consequence of Theorem~\ref{1}.

  \begin{theorem}\label{HBx}
  Let $h \in  {\cal Z}'_{+}$. Then
  \[
  N_{\pm}(h)=N_{\pm}(b)=N_{\pm}(s).
  \]
   \end{theorem}
  
  In particular, we have 
  
    \begin{theorem}\label{HBx1}
Let $h \in  {\cal Z}'_{+}$. Then the form $\pm {\pmb\la}h,\bar{f} \star f {\pmb\ra} \geq 0$ for all
 $f\in {\cal D}$ if and only if the form $\pm{ \la} b, {\cal J}\bar{g} * g{\ra}\geq 0$ for all
 $g\in C_{0}^\infty({\Bbb R})$, or the form $\pm{ \la} s , |u|^2{\ra} \geq 0$ for all $u \in {\cal Z} $. 
     \end{theorem}
     
 \medskip    
 
 {\bf 4.2.}
      In many cases the following consequence of Theorem~\ref{HBx} is convenient. According to Proposition~\ref{extb}, under the assumptions of Theorem~\ref{HBy}, $H$ is defined as the bounded self-adjoint operator corresponding to the form ${\pmb\la}h,\bar{f} \star f {\pmb\ra}$. Therefore $N_{\pm} (h)= N_{\pm} (H)$ is the total multiplicity of the  (strictly) positive spectrum  for the sign $``+$" and  of the (strictly) negative spectrum for the sign $``-$" of  the operator $H$.  For definiteness, we consider the negative spectrum.

       \begin{theorem}\label{HBy}
  Let $h \in  {\cal Z}_{+}'$, and let the corresponding sign-function satisfy condition \e{eq:LLs}.    If $s(x)\geq 0$, then the operator $H$ is positive. If  
  $s(x)\leq -s_{0}<0$  for almost all 
  $x $ in some interval $\Delta\subset{\Bbb R}$, then the operator $H$ has infinite negative spectrum.
   \end{theorem}
   
   \begin{pf}
   If $s(x)\geq 0$, then $H\geq 0$ according to the second relation \e{eq:M15}.
    
      Let $s(x)\leq -s_{0}<0$  for   $x \in \Delta$.  For an arbitrary $N$, we shall construct a linear set $ {\cal L}\subset  {\cal Z} $ of   dimension $N$ such that
   $s[u,u]<0$ for all $u\in  {\cal L}$, $u\neq 0$.  Then the second statement will follow from Theorem~\ref{HBx}.
   
   Choose a function $\varphi\in C_{0}^\infty ({\Bbb R})$ such that $\varphi(x)=1$ for $x\in [-\d,\d]$ and $\varphi(x)=0$ for $x\not\in [-2\d,2\d]$ where $\d=\d_{N}$ is a sufficiently small number. Let points $\alpha_{j}\in \Delta$, $j=1,\ldots, N$, be such that
   $\alpha_{j+1}-\alpha_{j} =\alpha_{j }-\alpha_{j-1}$ for $j=2,\ldots, N -1$. Set $\Delta_{j}= (\alpha_{j}-\d, \alpha_{j} +\d)$, $\wt{\Delta}_{j}= (\alpha_{j}-2\d, \alpha_{j} +2\d)$. For a sufficiently small $\d$, we may suppose that $\wt{\Delta}_{j}\subset \Delta$ for all $j=1,\ldots, N$ and that $\wt{\Delta}_{j+1}\cap \wt{\Delta}_{j}=\varnothing$ for   $j=1,\ldots, N -1$. We set $\varphi_{j}(x)=\varphi (x-\alpha_{j})$. Since $s(x)\leq -s_{0}<0$  for   $x \in \Delta$, we have   
    \begin{equation}
{ \la} s , |\varphi_{j} |^2{\ra} = \int_{-\infty}^\infty s (x) |\varphi_{j}(x) |^2 dx\leq -2\d s_{0}<0.
\label{eq:appa}\end{equation}
 The functions $\varphi_{1},\ldots,\varphi_N$ have disjoint supports and hence  ${ \la} s , |u |^2{\ra} < 0$ for an arbitrary non-trivial  linear combination $u$ of the functions $\varphi_{j}$.

The problem is that $\varphi_{j}\in C_{0}^\infty ({\Bbb R}) $ but $\varphi_{j}\not\in{\cal Z}$ (and even  $\varphi_{j}\not\in{\cal W}$). Thus we have to approximate these functions in the topology of $\cal S$ by functions $\varphi_{j}^{(\varepsilon)}\in{\cal Z}$. For a given $ \kappa >  1/2$ and an arbitrary $\varepsilon>0$, we can find a function  $\varphi^{(\varepsilon)}\in{\cal Z}$ such that
 \begin{equation}
\max_{x\in {\Bbb R}} \big((1+|x|)^{\kappa } |\varphi (x)-\varphi^{(\varepsilon)}(x)|\big) <\varepsilon.
\label{eq:app}\end{equation}
Now we set $\varphi_{j}^{(\varepsilon)} (x)=\varphi^{(\varepsilon)} (x-\alpha_{j})$. 

Let us check that for $\varepsilon$ small enough, the functions $\varphi_1^{(\varepsilon)},\ldots, \varphi_N^{(\varepsilon)}$ are linearly independent. Assume that
 \begin{equation}
 \sum_{j=1}^N  \lambda_{j} \varphi_{j}^{(\varepsilon)} (x)=0.
\label{eq:app1}\end{equation}
If $x\in \Delta_{k}$, then $|\varphi_k^{(\varepsilon)} (x)-1| <\varepsilon$ and
$|\varphi_j^{(\varepsilon)} (x) | <\varepsilon$ for $j\neq k$. Therefore \e{eq:app1} yields the estimate
 \[
(1-\varepsilon)|\lambda_{k}| \leq \varepsilon \sum_{j\neq k}  | \lambda_{j} |.
\]
Summing these estimates over $k=1,\ldots, N$, we see that
\[
(1-\varepsilon) \sum_{j=1}^N   | \lambda_{j} |  \leq \varepsilon (N-1)
 \sum_{j=1}^N | \lambda_{j} | .
 \]
 Hence $\lambda_{j}=0$ for all $j=1,\ldots, N$ if $\varepsilon (N -1)<1-\varepsilon$.

It follows from \e{eq:appa}, \e{eq:app} that, for all $j=1,\ldots, N$, 
\[
\int_{-\infty}^\infty s(x) |\varphi_{j}^{(\varepsilon)} (x)|^2 dx\leq -2\d s_{0}+ C \varepsilon
\]
and that 
\[
\big| \lambda_{j} \bar{\lambda}_k \int_{-\infty}^\infty s(x) \varphi_{j}^{(\varepsilon)} (x) \overline{\varphi_k^{(\varepsilon)} (x) } dx \big|\leq C \varepsilon (|\lambda_{j}|^2+|\lambda_k|^2) ,\q j\neq k.
\]
Thus we have the estimate
 \begin{equation}
\int_{-\infty}^\infty s(x) |\sum_{j=1}^N  \lambda_{j} \varphi_{j}^{(\varepsilon)} (x)|^2 dx \leq
- \sum_{j=1}^N |\lambda_{j}|^2 (2\d s_{0}- C \varepsilon).
\label{eq:app2}\end{equation}
The right-hand side here is negative if $\varepsilon$ is small enough and $\sum_{j=1}^N |\lambda_{j}|^2\neq 0$.

Let $ {\cal L}$ be spanned by the functions $\varphi_{1}^{(\varepsilon)},\ldots, \varphi_{N}^{(\varepsilon)}$ for sufficiently small $\varepsilon$. Then $\dim {\cal L}= N$ and $s[u,u]<0$ for all $u\in  {\cal L}$, $u \neq 0$ according to \e{eq:app2}. It follows from Theorem~\ref{HBx} that $N_{-}(H)\geq N$, and hence $N_{-}(H) =\infty$.
    \end{pf}

    Theorem~\ref{HBy} can be reformulated, although in a weaker   form, in terms of the functions $b(\xi)$ and even $h(t)$.   Suppose, for example,  that 
\begin{equation}
 b  \in L^1 ({\Bbb R}).
 \label{eq:N3}\end{equation}
 Then its Fourier transform $s(x)$ is a continuous function which tends to $0$ as $|x |\to \infty$. The operator $B$ defined by formula \e{eq:BB}  is  bounded, self-adjoint and
 \[
 \spec(B)=[\min_{x\in{\Bbb R}} s(x), \max_{x\in{\Bbb R}}s(x)].
 \]
 The  result below follows directly from Theorem~\ref{HBy}. Note that by Proposition~\ref{extb}  under assumption \e{eq:N3} the operator $H$ is compact.
  
  \begin{proposition}\label{VW}
  Under assumption \e{eq:N3}
 the Hankel operator $H$ is positive if and only if 
 $ s(x)\geq 0$. If $\min_{x\in{\Bbb R}} s(x)<0$, then necessarily $H$ has an infinite negative spectrum.
 \end{proposition}

  In particular, condition \e{eq:N3} is satisfied if
  \[
 h(t)=\frac{\theta(\ln t)}{t}\q {\rm where}\q \theta\in {\cal Z} .
\]
 In this case $a=\Phi\theta \in C_{0}^\infty ({\Bbb R})$ and hence $b\in C_{0}^\infty ({\Bbb R})$.   
 
 \medskip    
 
 {\bf 4.3.}
 For a proof that a Hankel operator is not sign-definite it is sometimes even not necessary to calculate the sign-function $s(x) $ (the Fourier transform   of   $b(\xi)$). It turns out that if
 $b(\xi)$      grows as $|\xi |\to \infty$, then      the  form $b[g,g]={ \la} b, {\cal J}\bar{g} * g{\ra}$ cannot be sign-definite. More precisely, we have the following statement about  convolutions   with growing kernels $b(-\xi)=\overline{b(\xi)}$.

\begin{theorem}\label{neq}
Let   $b=b_{0} + b_{\infty}$ where $b_{0}\in C^p ({\Bbb R})' $ for some $p \in {\Bbb Z}_{+}$ and $b_{\infty}\in L^\infty_{\rm loc}({\Bbb R}) $.
  Suppose that there exists a sequence of intervals $\Delta_{n}=(r_{n}-\sigma_{n}, r_{n}+\sigma_{n})$ where $r_{n}\to\infty$ $($or equivalently $r_{n}\to-\infty)$ and the sequence $\sigma_{n}$ is bounded   such that    
 \begin{equation}
\lim_{n\to \infty} \sigma_{n}^l \min_{\xi\in\Delta_{n}} \Re b_{\infty}(\xi)=  \infty
\q {\rm or}\q \lim_{n\to \infty} \sigma_{n}^l \max_{\xi\in\Delta_{n}} \Re b_{\infty}(\xi)= -\infty,
\label{eq:neq}\end{equation}
where $l=2$ if $p=0$ or $p =1$ and $l=p+1$ if $p \geq 2$. 
Then for both signs $N_{\pm} (b) \geq 1$.
 \end{theorem}  
 
 \begin{pf} 
 Since     $b$ can be replaced by $-b$, we can  assume that, for example,   the first  condition \e{eq:neq} is satisfied.
 Pick a real even function $\varphi\in C_{0}^\infty ({\Bbb R})$ such that $\varphi (\xi)\geq 0$, $\varphi  (\xi)=1$ for $|\xi| \leq 1/4$, 
$\varphi (\xi)=0$ for $|\xi | \geq 1/2$  and set
 \begin{equation}
g_{ n}(\xi)=\varphi  (( \xi - r_{n}/2)/\sigma_{n} )\pm \varphi(( \xi + r_{n}/2)/\sigma_{n} ).
\label{eq:neq2t}\end{equation}
An easy calculation shows that
 \begin{equation}
(({\cal J}g_{ n})* g_{n})(\xi)= 2 \sigma_{n} \psi (\xi/ \sigma_{n}) \pm \sigma_{n}\psi ((\xi-r_{n})/ \sigma_{n}) \pm \sigma_{n} \psi ((\xi+ r_{n})/ \sigma_{n})
\label{eq:neq2}\end{equation}
where $\psi= ({\cal J}\varphi ) * \varphi \in C_{0}^\infty ({\Bbb R})$. The function $\psi(\xi)$ is also even, $\psi (\xi)\geq 0$, $\psi  (\xi)\geq 1/8$ for $|\xi | \leq 1/8$ and
$\psi (\xi)=0$ for $|\xi | \geq 1$.

Since $|\la b_{0}, g\ra | \leq C \| g\|_{C^p}$, it follows from \e{eq:neq2} that  
 \begin{equation}
|\la b_{0}, ({\cal J}g_{ n})* g_{n} \ra | \leq  C \sigma_{n}^{1-p}. 
\label{eq:neq2z}\end{equation}
Moreover, again according to \e{eq:neq2} we have
 \begin{equation}
\la b_\infty, ({\cal J}g_{ n})* g_{n} \ra =2\sigma_{n}^2 \int_{-\infty}^\infty b_{\infty}(\sigma_{n}\eta) \psi(\eta )d \eta
\pm 2   \sigma_{n}^2 \int_{-\infty}^\infty \Re b_{\infty}(\sigma_{n}\eta + r_{n}) \psi(\eta)d \eta. 
\label{eq:neq2z1}\end{equation}
The first term in the right-hand side is $O(\sigma_{n}^2)$. For the second one, we use the estimate
 \begin{equation}
32    \int_{-\infty}^\infty \Re b_{\infty}(\sigma_{n} \eta +r_{n}) \psi(\eta)d \eta\geq  \min_{|\xi-r_{n}| \leq \sigma_{n} } \Re b_{\infty}(\xi).
\label{eq:neq1}\end{equation}

Let us first choose the sign $``+ "$ in \e{eq:neq2t}. Then using representation \e{eq:neq2z1} and putting together estimates \e{eq:neq2z} and \e{eq:neq1},  we obtain the lower bound
 \[
\la b , ({\cal J}g_{ n})* g_{n} \ra \geq -c (\sigma_{n}^{1-p} + \sigma_{n}^2)+\sigma_{n}^2/ 16 \min_{|\xi-r_{n}| \leq \sigma_{n}}   \Re b_{\infty}(\xi).
\]
If $p=0$ or $p=1$, then 
under the first  condition \e{eq:neq} the right-hand side here tends to $+\infty$ as $n\to\infty$. If $p \geq 2$, it is bounded from below by
\[
\sigma_{n}^{1-p} \big(-c+ \sigma_{n}^l/ 16 \min_{|\xi-r_{n}| \leq \sigma_{n}}   \Re b_{\infty}(\xi) \big)
\]
where the expression in the brackets tends again to $+\infty$ . Therefore $\la b , ({\cal J}g_{ n})* g_{n} \ra > 0$ for sufficiently large $n$.
 Similarly choosing the sign $``- "$ in \e{eq:neq2t}, we see that $\la b , ({\cal J}g_{ n})* g_{n} \ra <0$ for sufficiently large $n$.
 \end{pf}
 
 \begin{corollary}\label{neqx}
 Instead of condition \e{eq:neq} assume that
  \[
\lim_{|\xi| \to \infty}  \Re b_{\infty}(\xi)=  \infty
\q {\rm or}\q\lim_{|\xi| \to \infty}  \Re b_{\infty}(\xi) = -\infty.
\]
Then for both signs $N_{\pm} (b) \geq 1$.
 \end{corollary}

In contrast to  Theorem~\ref{neq} there are no restrictions in Corollary~\ref{neqx}  on the parameter $p$ in the assumption  $b_{0}\in C^p({\Bbb R})' $.

On the other hand,
  condition \e{eq:neq} permits   $\Re b(\xi) $ to tend to $\pm\infty$ only on some system of intervals. Moreover, the lengths of these intervals may tend to zero. In this case, however, the growth of   $\Re b(\xi) $ and the decay of these lengths should be correlated and there are restrictions on admissible values of the parameters $p$ and $l$.

Unlike Theorem~\ref{HBy}, Theorem~\ref{neq} does not guarantee that $N=\infty$; see subs.~5.4, for a discussion of various possible cases.

\medskip

 {\bf 4.4.}
  Theorem~\ref{HBx1}  can be combined with the Bochner-Schwartz theorem (see, e.g., Theorem~3 in \S 3 of Chapter II of the book \cite{GUEVI}). It states that a distribution $b\in C_{0}^\infty({\Bbb R})'$ satisfying the condition  ${ \la} b, {\cal J}\bar{g} * g{\ra} \geq 0$ for all
 $g\in C_{0}^\infty({\Bbb R})$  (such   $b$ are sometimes called    distributions     of positive type)  is the Fourier transform  
    \[
b(\xi) =(2\pi)^{-1 } \int_{-\infty}^\infty e^{-ix\xi} d{\sf M}(x )
\]
of a positive measure $d{\sf M}(x )$ such that
  \begin{equation}
\int_{-\infty}^\infty (1+|x|)^{-\varkappa } d{\sf M}(x )<\infty
\label{eq:Sch}\end{equation}
for some $\varkappa $ (that is,  of at most polynomial growth at infinity). In particular, this ensures that $b\in{\cal S}'$.

Theorem~\ref{HBx1} implies that if ${\pmb\la}h,\bar{f} \star f {\pmb\ra} \geq 0$ for all
 $f\in {\cal D}$, then the corresponding  distribution $b$ is        of positive type. It means that the  sign-function $s(x)$ of $h(t)$ is determined by the measure $d{\sf M}(x )$:
  \[
\la s, \varphi\ra =  \int_{-\infty}^\infty \ov{\varphi(x)} d{\sf M}(x ), \q \varphi \in {\cal S},
\]
that is, $s(x)dx=d{\sf M}(x )$.  Let us define the measure
  \begin{equation}
  d m(\lambda)= \lambda d{\sf M}(-\ln\lambda),\q \lambda\in{\Bbb R}_{+}.
 \label{eq:Mm}\end{equation}
  It is easy to see that condition \e{eq:Sch} is equivalent to   condition \e{eq:Sch1}
on    measure \e{eq:Mm}. In terms of distribution   \e{eq:B}, we have $\lambda h^\natural (\lambda)d\lambda= dm(\lambda)$. 
 Therefore Theorem~\ref{round} leads to the following result.

\begin{theorem}\label{HBx2}
Let $h \in {\cal Z}_{+}'$ and ${\pmb\la}h,\bar{f} \star f {\pmb\ra}\geq 0$ for all
 $f\in {\cal D}$. Then $h(t)$ admits the representation \e{eq:Conv} with
 a   positive measure $dm (\lambda )$, $\lambda\in {\Bbb R}_{+}$, satisfying for some $\varkappa $    condition \e{eq:Sch1}. 
      \end{theorem}

      The representation \e{eq:Conv} is of course  a particular case of \e{eq:conv1}. It is much more precise than \e{eq:conv1} but requires the positivity of ${\pmb\la}h,\bar{f} \star f {\pmb\ra}$. Theorem~\ref{HBx2} shows that  the positivity of ${\pmb\la}h,\bar{f} \star f {\pmb\ra}$ imposes very strong conditions on $h(t)$. Actually, we have
      
      \begin{corollary}\label{HBx3}
       Let $h \in  {\cal Z}'_{+}$ and ${\pmb\la}h,\bar{f} \star f {\pmb\ra}\geq 0$ for all
 $f\in {\cal D}$.  Then $h\in C^\infty ({\Bbb R}_{+})$ and 
    \begin{equation}
(-1)^n h^{(n)}(t)\geq 0
\label{eq:CoMo}\end{equation}
for all $t>0$ and all $n=0,1,2,\ldots$
 $($such functions are   called completely monotonic$)$. The function  $h(t)$ admits an analytic continuation in the right-half plane $\Re t> 0$ and it is uniformly bounded in every strip $\Re t\in (t_{1},t_{2})$ where $0< t_{1} <t_{2} <\infty$. Moreover, for some $\varkappa\in {\Bbb R}$ and $C>0$ the estimate holds:
    \begin{equation}
 h(t) \leq C t^{-1} (1+|\ln t|)^\varkappa ,\q t>0.
\label{eq:CoMo1}\end{equation}
      \end{corollary} 
      
      All these assertions are direct consequences of the representation \e{eq:Conv}. In particular, under condition  \e{eq:Sch1} we have
      \[
      h(t) \leq C \max_{\lambda\geq 0} \big(e^{-t\lambda}\lambda (1+|\ln \lambda|)^\varkappa\big)
      \]
      which yields \e{eq:CoMo1}.

 Note that according to the  Bernstein theorem (see, e.g., Theorems~5.5.1 and 5.5.2 in \cite{AKH}) condition \e{eq:CoMo} implies that the function $h(t)$ admits the representation \e{eq:Conv} with some measure $dm(\lambda)$. Of course, condition \e{eq:CoMo}  does not impose any restrictions on the measure $dm(\lambda)$ (except that the integral \e{eq:Conv} is convergent for all $t>0$). In contrast to the  Bernstein theorem we deduce the representation \e{eq:Conv} from the positivity of the Hankel form. In this context  condition  \e{eq:Sch1} is due to the assumption $h \in {\cal Z}'_{+}$.

  We mention also a related result of H.~Widom. He considered in \cite{Widom} Hankel operators $H$ with kernels $h(t)$ admitting the representation \e{eq:Conv} and showed that $H$ is bounded if and only if $m([0,\lambda))=O(\lambda)$ as $\lambda\to 0$ and as $\lambda\to \infty$.  In this case $h(t)\leq C t^{-1}$ for some $C>0$. Thus Theorem~\ref{HBx2} and Corollary~\ref{HBx3} can be regarded as an extension of Widoms's results to unbounded operators.
   
   Under the positivity assumption the identity \e{eq:M15f} takes a more precise form.  
  
   \begin{proposition}\label{1d}
Let $h \in  {\cal Z}'_{+}$ and   ${\pmb\la}h,\bar{f} \star f {\pmb\ra} \geq 0$ for all
 $f\in {\cal D}$. Then there exists  a  positive measure $d{\sf M}(x )$ satisfying condition  \e{eq:Sch} for some $\varkappa $ such that 
    \[
{\pmb\la}h,\bar{f}_{1} \star f_{2} {\pmb\ra} =  \int_{-\infty}^\infty u_{1}(x) \overline{u_{2}(x)} d{\sf M}(x )
\]
for all $f_{j}\in {\cal D}$, $j=1,2$, and $u_{j}=\Phi^* \Xi f_{j}\in {\cal Z}$ where the mapping $\Xi$ is defined by \e{eq:MAID1}.
 \end{proposition}

  %%%%%%%%%%%%%%%%%%%%%%%%%%%%%%
 %************************************************************
\section{Applications and examples}  
%**%%%%%%%%%%%%%%%%%%%%%%%%%%%

 ********************************************************

 {\bf 5.1.}
 Consider first  self-adjoint Hankel operators $H$ of  finite rank. Recall that
integral kernels of   Hankel operators   of  finite rank are given (this is the Kronecker theorem -- see, e.g., Sections~1.3 and 1.8 of the book \cite{Pe}) by the formula
\begin{equation}
h(t) =\sum_{m=1}^M  P_{m} (t)e^{- \alpha_{m}  t}
\label{eq:FDvm}\end{equation}
where $\Re\alpha_{m}  >0$ and 
 $   P_{m} (t)$ 
    are polynomials of degree $K_{m} $. If $H$ is self-adjoint, that is, $ h(t)=\ov{h(t)} $,    then the set $\{\alpha_1,\ldots,\alpha_{M}\}$ consists of points lying on the real axis and pairs  of points symmetric with respect to it.   Let $\Im \alpha_{m}=0$ for $m=1,\ldots, M_{0}$ and  $\Im \alpha_{m}>0$,  $\alpha_{ M_{1}+m}=\bar{\alpha}_m$ for $m=M_{0}+ 1,\ldots, M_{0}+M_{1}$. Thus $M=M_{0}+2M_{1}$; of course the cases $M_{0}=0$ or $M_{1}=0$ are not excluded. The condition  $h(t)=\ov{h(t)} $ also requires that $ P_{m} (t) =\ov{P_{m} (t)}  $ for $m=1,\ldots, M_{0}$ and $ P_{M_{1}+m}  (t)= \ov{P_{ m} (t)} $  for $m=M_{0} +1, \ldots,M_{0} + M_{1} $.    As is well known and as we shall see below,
      \[
\rank H=\sum_{m=1}^M K_{m}+M=: r.
\]
For $m= 1, \ldots,M_{0}  $, we denote by
      ${\sf p}_{m}=\bar{\sf p}_m$   the coefficient at $t^{K_{m}}$ in the polynomial
       $P_{m}(t)$.       

  The following assertion yields an explicit formula for the numbers $N_\pm (H)$. Its proof will be given in \cite{Y3}.
 
  \begin{theorem}\label{FDH1}
For $m= 1, \ldots,M_{0}  $, set
  \begin{equation}
\left. \begin{aligned}
{\cal N}_{+}^{(m)}  = {\cal N}_{-}^{(m)} &=(K_{m} +1)/2 \q {\rm if} \; K_{m} \; {\rm is  } \; {\rm   odd}
 \\
{\cal N}_{+}^{(m)}  -1= {\cal N}_{-}^{(m)}  &= K_{m}  /2 \q {\rm if} \; K_{m} \; {\rm is  } \; {\rm   even} \; {\rm   and}
 \q {\sf p}_{m}  > 0
 \\
{\cal N}_{+}^{(m)}= {\cal N}_{-}^{(m)}  -1  &= K_{m}  /2  \q {\rm if} \; K_{m} \; {\rm is  } \; {\rm   even} \; {\rm   and}
 \q {\sf p}_{m} < 0 .
    \end{aligned}
    \right\}
    \label{eq:RXm}\end{equation} 
  Then the total numbers $N_{\pm} (H)$  of   $($strictly$)$ positive and negative eigenvalues of the operator $ H$ are given by the formula 
       \begin{equation}
 N_{\pm} (H) =  \sum_{m=1}^{M_{0}}  {\cal N}_\pm^{(m)} + \sum_{m=M_{0}+1}^{M_{0}+M_1}  K_{m} + M_{1} .
   \label{eq:TN}\end{equation}
 \end{theorem}

 Formula \e{eq:RXm} shows that every pair 
      \begin{equation}
 P_{m} (t) e^{-\alpha_{m}t} + P_{m+M_{1}} (t) e^{-\alpha_{m+M_{1}}t},\q m= M_{0}+1,\ldots, M_{0}+M_{1},
   \label{eq:FDvr}\end{equation}
   of complex conjugate terms in \e{eq:FDvm}   yields $K_{m}+1$ positive and $K_{m}+1$ negative eigenvalues. The contribution of every real term $P_{m} (t)e^{- \alpha_{m}  t}$ where $m= 1,\ldots, M_{0} $ also consists of the equal numbers $(K_{m}+1)/2$ of positive and  negative eigenvalues if the degree $K_{m}$ of the polynomial $P_m (t)$ is odd. If $K_{m}$   is even, then there is one more positive (negative) eigenvalue if ${\sf p}_{m}>0$
 (${\sf p}_{m}<0$).  In particular,  in the question considered, there is no ``interference" between different terms $P_{m} (t) e^{-\alpha_{m}t}$, $m=1,\ldots, M_{0}$, and pairs \e{eq:FDvr} in representation \e{eq:FDvm} of the kernel $h(t)$.

According to \e{eq:TN} the operator $H$ cannot be sign-definite if $M_{1}>0$.  Moreover, according to \e{eq:RXm}, ${\cal N}_{\pm}^{(m)}=0$ for $m=1,\ldots, M_{0}$ if and only if $K_{m}=0$ and $\mp {\sf p}_{m}>0$. Therefore we have the following result.
   
  \begin{corollary}\label{FD}
A  Hankel operator  $H$ of finite rank   is positive  $($negative$)$ if
and only if its kernel is given by the formula 
   \[
h(t) =\sum_{m= 1}^{M_{0}} {\sf p}_{m} e^{- \alpha_{m} t}
\] 
where $\alpha_{m}>0$ and $  {\sf p}_{m}>0$ $({\sf p}_{m}<0)$.
 \end{corollary}
 
 Corollary~\ref{FD} admits different proofs which avoid formula \e{eq:TN}. For example, one can use that although the functions $P_{m}(t)e^{-\alpha_{m}t}$ are analytic   in the right-half plane $\Re t>0$,   they are bounded for $t=\tau+i\sigma$ as $\sigma\to \infty$ for a constant $P_{m}(t)$ only. Therefore according to Corollary~\ref{HBx3}  such Hankel operators cannot be positive. Alternatively, using formula \e{eq:E33} below,  one can  deduce Corollary~\ref{FD}    from Theorem~\ref{neq}.  
     
\medskip
 
{\bf 5.2.}
Consider now   Hankel operators   $H$ with kernels \e{eq:E1r}. Since  the case $k=0,1,\ldots$ (finite rank Hankel operators) has been  discussed   in the previous subsection, here we suppose that $k\neq 0,1,\ldots$.
If $k > -1$,  condition \e{eq:ass} is   satisfied for all  $\kappa $, and the operators   $H$ are compact (actually, they belong to much better classes of operators).   If $k=-1$, then condition \e{eq:ass}  is   satisfied for    $\kappa >1$, and the operators   $H$ are bounded but not compact. 

Let us calculate the corresponding $b$- and $s$-functions.  
 If $k > -1$, then    function \e{eq:M6}  equals
 \begin{equation}
a(\xi) = (2\pi)^{-1/2} \int_{0}^\infty t^k e^{-\alpha t} t^{-i\xi} dt =
(2\pi)^{-1/2} \alpha^{-1-k + i \xi} \Gamma (1+k-i \xi),
\label{eq:E2}\end{equation}
 and hence  function \e{eq:M9}  equals
 \begin{equation}
b(\xi) =  \alpha^{-1-k + i \xi} \frac{\Gamma (1+k- i \xi)} { 2\pi \Gamma (1 -i \xi)}.
\label{eq:E3}\end{equation}
 
If  $k=-1$, then in accordance with formulas \e{eq:E2} and \e{eq:E3}, we have
\[
a(\xi) = (2\pi)^{-1/2} \alpha^{ i \xi} \lim_{\varepsilon\to +0}\Gamma (\varepsilon-i \xi),\q
b(\xi) = (2\pi)^{-1 } \alpha^{   i \xi}i (\xi+i0)^{-1}.
\]
It yields the expression
 \begin{equation}
s (x) = 0, \q x > \beta, \q s (x) = 1, \q x < \beta,\q{\rm where}\q \beta=- \ln \alpha,
\label{eq:E3c}\end{equation}
for the function $s=\sqrt{2\pi}\Phi^* b$.  Formula \e{eq:E3c} remains true for the Carleman operator $\bf C$ (the Hankel operator with kernel $h(t)=t^{-1}$) when $\alpha=0$. Indeed, in this case  according to \e{eq:carl1}   the sign-function $s(x)=1$.   
       
Next, we calculate the Fourier transform   of function \e{eq:E3}.  
    Assume first that $k\in (-1,0)$. Then (see, e.g.,  formula (1.5.12) in the book \cite{BE})
   \[
   \int_{0}^\infty t^{-k-1} (t+1)^{-1+i\xi} dt =\frac{\Gamma(-k)\Gamma(1+k-i\xi)}{ \Gamma(1-i\xi)}.
   \]
   Making here the change of variables $t+1= \alpha^{-1} e^{-x}$, we find that
     \[
\frac{1}{\Gamma(-k)}   \int_{-\infty}^\infty (e^{-x} -\alpha)_{+}^{-k-1} e^{-ix\xi} dx =
\alpha^{-1-k-i\xi} \frac{ \Gamma(1+k-i\xi)}{ \Gamma(1-i\xi)}.
   \]
  Passing now to the inverse Fourier transform, we see  that for $k\in (-1,0)$ the sign-function $s(x)= s_{k}(x)$ of kernel \e{eq:E1r}  equals
   \begin{equation}
   s (x)=
   \frac{1}{\Gamma(-k)}   (e^{-x} -\alpha)_{+}^{-k-1} .
\label{eq:bb}\end{equation}

Let us verify that this formula remains true for all non-integer $k$. To that end, we assume that \e{eq:bb} holds for some non-integer $k> -1$ and check   it   for $k_{1}=k+1$. Since
\[
\Gamma(1+k_{1}-i\xi)=(k_{1}-i\xi)\Gamma(1+k-i\xi),
\]
we have
\[
s_{k_{1}}(x)=\alpha^{-1} (k_{1}-\partial)s_{k}(x).
 \]
 Substituting here formula \e{eq:bb} for $s_{k }(x)$ and differentiating this expression, we obtain formula
 \e{eq:bb} for $s_{k_{1}}(x)$. This concludes the proof of relation
 \e{eq:bb}   for all $k\geq -1$.

 \begin{lemma}\label{Hks}
Let $h(t)$ be given by formula \e{eq:E1r} where $ k\not\in {\Bbb Z}_{+}$.  Then the sign-function is determined by relation \e{eq:bb}.
 \end{lemma}
 
 Actually,   relation  \e{eq:bb} remains true for $k\in {\Bbb Z}_{+}$ if one takes into account that the distribution $  (e^{-x} -\alpha)_{+}^{-k-1}$ has poles at integer points. 
 For example, for $k=0$ we have
 \begin{equation}
   s (x)=
\alpha^{-1} \d (x-\beta).
\label{eq:bbrr}\end{equation}

Obviously, $s (x)=0$ for $x>\beta=-\ln\alpha$. If $k=-1$, then $s (x)=1$ for $x<\beta $. If $k\in (-1,0)$, then   $s (x)\geq 0$ and  $s\in L^1({\Bbb R})$. Therefore it follows from Theorem~\ref{HBy} that $H\geq 0$.
 
 If $k>0$, then  distribution  \e{eq:bb} does not have a definite sign. Therefore it can be deduced from Theorem~\ref{HBx1} that the corresponding   Hankel operator also  is not sign-definite. 
 
 Alternatively, for the proof of this result we can use Corollary~\ref{neqx}.   
Formula \e{eq:M11}   implies that function \e{eq:E3}    has the asymptotics 
   \begin{equation}
b(\xi)= (2\pi)^{-1}\alpha^{-1-k-i\xi} (-i \xi)^k (1+ O(|\xi|^{-1})),\q |\xi|\to \infty.
\label{eq:E8}\end{equation}
Making the dilation transformation in \e{eq:E1r}, we can suppose that $\alpha=1$.
Then we have
 \begin{equation}
\Re b(\xi)= (2\pi)^{-1} \cos (\pi k/2) \xi^k  + O(\xi^{k-1}),\q \xi\to +\infty,
\label{eq:E9}\end{equation}
Since $ \cos (\pi k/2) \neq 0$  unless $k$ is an integer odd number, 
this expression tends to $\pm\infty$ if $\pm\cos (\pi k/2) >0$.   Thus Corollary~\ref{neqx} for the case $b_{0}=0$  ensures that the Hankel operator $H$ is not sign-definite.

Let us summarize the results obtained.

 \begin{proposition}\label{HK}
  The Hankel operator with kernel \e{eq:E1r}  is positive for   $k\in [-1,0]$, and it is not sign-definite  for   $k>0$.
  \end{proposition}

Explicit formulas for the sign-functions can also be used to treat   more complicated Hankel operators. For example,
 in view of   \e{eq:E3c} the following assertion directly follows from Theorem~\ref{HBy}.

 \begin{example}\label{car}
 The Hankel operator with kernel
 \[
 h(t) =t^{-1} (e^{-\alpha_{1}t}-\gamma e^{-\alpha_{2}t}),\q \gamma\geq 0,
 \]
is positive if and only if $\alpha_{2}\geq \alpha_{1}\geq 0$ and $\gamma\leq 1$.
 \end{example}

 \medskip
 
 {\bf 5.3.}
In this subsection, we consider  the  Hankel operator $H$ with  kernel \e{eq:F5}.   Condition \e{eq:ass}  is now fulfilled for all $\kappa$, and the   operator $H$  belongs of course to the Hilbert-Schmidt class (actually, to much better classes).
Observe that
 \[
 a(\xi)=(2\pi)^{-1/2} \int_{0}^\infty e^{-t^r} t^{-i \xi} dt= (2\pi)^{-1/2} r^{-1} \Gamma( (1-i \xi)/r)
 \]
 and define, as usual, the function $ b(\xi)$   by formula \e{eq:M9} so that
 \begin{equation}
b(\xi)= (2\pi  r)^{-1}  \frac{\Gamma( (1-i \xi)/r)}{\Gamma( 1-i \xi) } .
\label{eq:Exp}\end{equation}
 
 Consider first the case $r>1$.  It follows from  the Stirling formula \e{eq:M11} that for all $r>1$ the modulus of   function \e{eq:Exp}
     exponentially grows and the periods of its oscillations tend to zero only logarithmically
 as $| \xi |\to\infty$. Therefore Theorem~\ref{neq}   implies that the Hankel operator with kernel  \e{eq:F5} is not  sign-definite.  
   
  The Hankel operator $H$ with kernel
$h(t)=e^{-t^2}$ can also be treated (see Appendix~B) in a completely different way which is perhaps  also of some interest. This method shows that both positive and negative spectra of the  operator $H$ are infinite.

If $r=1$, then $h(t)=e^{-t}$ yields a positive Hankel operator of rank $1$.  

Let us now consider the case $r<1$. Then, again according to the Stirling formula \e{eq:M11},    function \e{eq:Exp} belongs to $L^1 ({\Bbb R})$ so that its Fourier transform  
  \begin{equation}
s(x)=(2\pi r)^{-1} \int_{- \infty}^\infty \frac{\Gamma( (1-i \xi)/r)}{\Gamma( 1-i \xi) }e^{ix\xi}d\xi =: I_r (x)
\label{eq:exp}\end{equation}
 is a continuous function which tends to $0$ as $|x|\to \infty$. Therefore by Proposition~\ref{VW}the corresponding Hankel operator $H\geq 0$ if and only if $I_r (x)\geq 0$ for all $x\in{\Bbb R}$.
 
 It turns out that
   $I_{r}(x)\geq 0$. Surprisingly, we have not found a proof of this fact in the literature, but it follows from our results. Only for $r=1/2$, integral \e{eq:exp} can be explicitly calculated. Indeed, 
  according to formula (1.2.15) of \cite{BE}
 \[
   \frac{ \Gamma(2 (1-i\xi))}{ \Gamma( 1-i\xi )}=  2^{1- 2 i \xi}  \pi^{-1/2}  \Gamma( 3/2-i \xi) .
 \]
 Therefore it follows from formula \e{eq:IDGg}  that
 \begin{equation}
 I_{1/2} (x)= 2^{-1} \pi^{-1/2}  e^{3x/2} e^{-e^x/4}
 \label{eq:expZ}\end{equation}
 which is of course positive.
 
 For an arbitrary $r\in (0,1)$, one can proceed from the Bernstein theorem on completely monotonic functions (see subs.~4.4). Observe that if 
 \begin{equation}
 \psi(t)=t^{-p} e^{-t^r}, \q p\geq 0,
\label{eq:exp1}\end{equation}
 then
  \[
\psi'(t)=-p t^{-p-1} e^{-t^r} -r t^{-p+r-1} e^{-t^r}.
 \]
 Further differentiations of $\psi(t)$ change the sign and yield sums of terms having the form \e{eq:exp1}. Thus the function $h(t)=e^{-t^r}$ satisfies for all $n$   condition \e{eq:CoMo} and hence admits the representation \e{eq:Conv} with some positive measure $dm(\lambda)$. It follows from \e{eq:Conv} that 
 \[
(Hf,f)=   \int_{0}^\infty |({\sf L}f)(\lambda)|^2 dm(\lambda)\geq 0, \q \forall f\in C_{0}^\infty ({\Bbb R}_{+}),
 \]
where  ${\sf L}$ is the Laplace transform \e{eq:LAPj}. Since the operator $H$ is bounded, this implies that $H\geq 0$.
 
 Thus we have obtained the following result.
 
   \begin{proposition}\label{exp}
 The Hankel operator with kernel \e{eq:F5}  is positive for   $r\in (0,1]$, and it is not sign-definite  for  $r>1$.
 \end{proposition}
 
 Putting together this result with Theorem~\ref{HBy}, we see that integral \e{eq:exp} is positive for all $r\in (0,1)$. Our indirect proof of this fact looks curiously enough.

 \medskip

 {\bf 5.4.}
  Let us now  discuss  convolution operators with growing kernels $b(\xi)$. We emphasize that condition \e{eq:neq} does not guarantee that the numbers $N_{\pm}(b)$ are infinite. Indeed, consider the kernel $h(t)=t^k e^{-\alpha t}$ where   $k$ is a positive integer.   Formula  \e{eq:E3} shows that for $\Im \alpha=0$ the corresponding function  $b(\xi)$
    \begin{equation}
b(\xi) = (2\pi)^{-1} \alpha^{-1-k + i \xi}   (1 -i \xi)\cdots (k -i \xi)
\label{eq:E33}\end{equation}
 has a power asymptotics as  $|\xi|\to\infty$.  According to Theorem~\ref{FDH1}
  the positive and negative spectra of the   Hankel operator $H$ with  the kernel $h(t)$ are finite; for example, $H$ has exactly $(k+1)/2$ positive and negative eigenvalues if $k$ is odd. Moreover, if $\Im \alpha\neq 0$, then in view of \e{eq:E33} the function  $b(\xi)$ exponentially grows as   $\xi \to +\infty$ or $\xi \to -\infty$. Nevertheless the   Hankel operator $H$ with kernel $h(t)=t^k (e^{-\alpha t} + e^{-\bar{\alpha} t})$ has exactly $k+1$ positive and negative eigenvalues.
  
  On the other hand, for kernel  \e{eq:Exp} where $r=2$ we have $N_{\pm}(b)=\infty$. This follows from Theorem~\ref{HBx}  because, by Proposition~\ref{BB}, the Hankel operator with kernel $h(t)=e^{-t^2}$ has infinite number of positive and negative eigenvalues.
  
  A similar phenomenon occurs for Hankel operators with non-smooth kernels. This is discussed in the next section.

 \section{Hankel operators with non-smooth kernels}

 {\bf 6.1.}
  Let the symbol of the   Hankel operator $H$ be defined by the formula $\omega (\mu)=e^{i t_{0} \mu}$. Evidently, $\omega \in {\tt H}^\infty_{_-}$ if $t_{0}\leq 0$, and hence $H=0$ in this case. If $t_{0}>0$, then the integral kernel of $H$ equals $\d (t -t_{0})$ so that
\[
(H f)(t)= f(t_{0}-t).
\]
Condition \e{eq:CHb} is now satisfied for $N=0$ and $\kappa=0$.

The operator $H$  admits an explicit spectral analysis. Indeed, observe first  that $(H f)(t)= 0$ for $t> t_{0}$ and hence $L^2 (t_{0},\infty) \subset \Ker H$.  
 Since $H^2 f= f$ for $f\in L^2 (0,t_{0})$,  the restriction of $H$ on its invariant subspace $L^2 (0,t_{0})$  may have only $\pm 1$ as eigenvalues. Obviously, the eigenspace ${\cal H}_{\pm}$ of  $H$ corresponding to the eigenvalue $\pm 1$  consists of all functions $f(t)$ such that  $  f(t)= \pm f(t_{0}-t)$. Since 
 $$
 {\cal H}_{+}\oplus {\cal H}_{-}\oplus L^2 (t_{0},\infty) = L^2 ({\Bbb R}_{+}),
 $$
  the spectrum of $H$ consists of the eigenvalues $0, 1, -1$ of infinite multiplicity each.

    In this example, the   $b$-function equals
  \begin{equation}
b(\xi)=  \frac{t_{0}^{-i\xi}}{ 2\pi \Gamma( 1-i \xi) }. 
\label{eq:Exp1}\end{equation}
  Note that all functions \e{eq:Exp} where $r>1$,   \e{eq:Exp1}  as well as \e{eq:Exp2} below exponentially grow and  oscillate at infinity. In these cases the corresponding Hankel operators have infinite positive and negative spectra.

 \medskip

 {\bf 6.2.}
 It follows  from
   Corollary~\ref{HBx3} that a Hankel operator $H$ can be sign-definite only for kernels  $h\in C^\infty({\Bbb R}_{+})$.    Actually, if $h(t)$ or one of its derivatives  $h^{(l)}(t)$ has a jump discontinuity, then   $H$ has infinite number of both positive $\lambda_{n}^{(+)}$ and negative $- \lambda_{n}^{(-)}$, $n=1,2, \ldots,$ eigenvalues,  and we     can even calculate their asymptotics as $n\to\infty$.  Positive   (negative) eigenvalues are of course  enumerated in decreasing (increasing) order with multiplicities taken into account.
   
   Let us start with an explicit 
   
    \begin{example}\label{powerex}
Let $h(t)=(t_{0} -t)^{l}$ for some $l=0, 1,\ldots$ if  $t\leq t_{0}$ and $h(t)=0$ if $t >t_{0}$. Then
   \[
(H f)(t) =  \int_{0}^{t_{0}-t} (t_{0} -t-s)^{l} f(s) ds,\q t\in (0,t_{0}),
\]
% \label{eq:EXE}\end{equation}
 and $(H f)(t) = 0$ for $t\geq t_{0}$. For such $h(t)$, 
  the symbol equals 
\[
\omega (\mu)=  l! (i\mu)^{-l-1} \big(  e^{i \mu t_{0}} -\sum_{k=0}^l \frac{1}{k!}(i\mu t_{0})^k\big)
 \]
and the  $b$-function  equals 
\begin{equation}
b(\xi)=  \frac{ l! t_{0}^{l+1-i\xi}}{ 2\pi \Gamma(l+ 2-i \xi) }. 
 \label{eq:Exp2}\end{equation}
  \end{example} 
 
  Let us consider the spectral problem $Hf=\lambda f$, that is, 
    \begin{equation}
 \int_{0}^{t_{0}-t} (t_{0} -t-s)^{l} f(s) ds= \lambda f(t),\q t\in (0,t_{0}).
 \label{eq:EXEs}\end{equation}
 Differentiating this equation $k$ times,   we find that
   \begin{equation}
(-1)^k l(l-1) \cdots (l-k+1) \int_{0}^{t_{0}-t} (t_{0} -t-s)^{l-k} f(s) ds= \lambda f^{(k)}(t) 
 \label{eq:EXEs1}\end{equation}
 for $k=1,\ldots, l$.
  Differentiating the last equation where $ k=l$ once more, we see that
    \begin{equation}
  l !  f(t)  = \lambda (-1)^{l+1} f^{(l+1)}(t_{0} -t),\q t\in (0,t_{0}).
 \label{eq:EXEs2}\end{equation} 
 Setting  $t=t_{0}$ in \e{eq:EXEs} and  \e{eq:EXEs1}, we obtain the boundary conditions
   \begin{equation}
  f(t_{0} )  = f' (t_{0} )=\cdots  = f^{(l)} (t_{0} ) =0.
 \label{eq:EXEbc}\end{equation} 
 Conversely, if a function $f(t)$ satisfies equation \e{eq:EXEs2} and boundary conditions
 \e{eq:EXEbc}, it satisfies also equation \e{eq:EXEs}. This leads to the following intermediary result.
 
  \begin{lemma}\label{powerex1}
Let the   operator $A$ be defined  on the Sobolev class ${\sf H}^{l+1}  (0,t_{0})$ by the equation
  \begin{equation}
 (A f)(t)= (-1)^{l+1} f^{(l+1)}(t_{0}-t).
 \label{eq:AAq}\end{equation}
Considered with boundary conditions
 \e{eq:EXEbc},  it is self-adjoint   in the space $L^2 (0,t_{0})$,  and its eigenvalues $\alpha_{n}$ are linked to eigenvalues $\lambda_{n}$ of the operator $H$ by the equation $\alpha_{n} = l!  \lambda_{n}^{-1}$.
 \end{lemma}
 
  \medskip

 {\bf 6.3.}
 Clearly,  $A^2$ is a differential  operator and the asymptotics of its eigenvalues is described by the Weyl formula. However, to find the asymptotics of   eigenvalues of the operator $A$, we have to distinguish its positive and negative eigenvalues. For this reason, it is convenient to introduce an auxiliary operator $\wt{A}$ with symmetric (with respect to the point $0$) spectrum having the same asymptotics of eigenvalues as $A$.
 
 We define the operator $\wt{A}$  by the same formula \e{eq:AAq}
as $A$ but consider it on functions in  $ {\sf H}^{l+1}  (0,t_{0}/2)\oplus {\sf H}^{l+1} (t_{0}/2,t_{0})$  satisfying the boundary conditions
 \begin{equation}
  f^{(k)}(0 )  =  f^{(k)}(t_{0}/2-0 ), \q  f^{(k)}(t_{0}/2 +0 ) =    f^{(k)}(t_{0} ),  
 \label{eq:EXEbc1x}\end{equation}
where $k=0,\ldots, l$ for $l$ even and
  \begin{equation}
  f^{(k)}(0 )  =  f^{(k)}(t_{0}/2-0 ) =f^{(k)}(t_{0}/2+0 )=    f^{(k)}(t_{0} )=0,  
 \label{eq:EXEbc1y}\end{equation}
where $ k=0,\ldots, (l-1)/2$ for $l$ odd. The operator $\wt{A}$ is self-adjoint 
   in the space
 \[
  L^2 (0,t_{0})=L^2 (0,t_{0}/2)\oplus L^2 (t_{0}/2,t_{0}),
\]
% \label{eq:oplus}\end{equation}
and it is determined   by the matrix
  \begin{equation}
\wt{A}= \begin{pmatrix}
0 & A_{1,2}
\\
A_{2,1}    & 0
\end{pmatrix}  , \q  A_{1,2} = A_{2,1}^*,
 \label{eq:oplus1}\end{equation}
where $A_{2,1}: L^2 (0,t_{0}/2)\to  L^2 (t_{0}/2,t_{0}) $. The operator $A_{2,1}$ is  again given  by relation \e{eq:AAq} on functions in $ {\sf H}^{l+1} (0, t_{0}/2 )$ satisfying conditions \e{eq:EXEbc1x} or \e{eq:EXEbc1y} at the points $0$ and $t_{0}/2 -0$. It follows from formula \e{eq:oplus1} that the spectrum of the operator $\wt{A}$ is symmetric with respect to the point $0$ and consists of eigenvalues $\pm a_{n}$ where
  $a_{n}^2$  are eigenvalues of the operator $A_{2,1}^*A_{2,1}=:{\bf A}$. 

An easy calculation shows that    ${\bf A} $ is
 the differential operator ${\bf A}=(-1)^{l+1}  \partial^{2l+2}$   in the space $L^2 (0,t_{0}/2)$ defined  on functions in  the class  $ {\sf H}^{2l+2} (0,t_{0}/2)$
   satisfying the boundary conditions $ f^{(k)}(0 )  =  f^{(k)}(t_{0}/2 ) $ where $ k=0,\ldots, 2l+1$  for $l$ even and the boundary conditions 
 $$
  f^{(k)}(0 )  =  f^{(k)}(t_{0}/2 )=  f^{(l+1+k)}(0 )  =  f^{(l+1+ k)}(t_{0}/2 )= 0,
  $$
where $ k=0,\ldots, (l-1)/2$  for $l$ odd.  The asymptotics of   eigenvalues $a_{n}^2 $ of ${\bf A} $  is given by the Weyl formula, that is,
   \[
a_{n} =   (2\pi t_{0}^{-1} n)^{l+1} (1+O (n^{-1})).
\] 
 
 Let us now observe that the operators $A$ and $\wt{A}$ are self-adjoint extensions of a symmetric operator $A_{0}$ with finite deficiency indices. For example, $A_{0}$ can be defined by formula \e{eq:AAq} on $C^\infty$ functions vanishing in some neighbourhoods of the points $0$, $t_{0}/2$ and $t_{0}$.
  Therefore the operators $A$ and $\wt{A}$ have the same asymptotics of spectra. Taking  Lemma~\ref{powerex1} into account, we obtain the following result.
 
   \begin{lemma}\label{powerex2}
   Eigenvalues of the Hankel operator $H$ defined in
   Example~\ref{powerex} have the asymptotics 
      \[
\lambda^{(\pm)}_{n} = l!  (2\pi)^{-l-1}  t_{0}^{l+1} n^{-l-1} (1+O (n^{-1})).
\]
 \end{lemma}
   
  \begin{remark}\label{pwerex2}
   In the case $l=0$ we have the explicit formulas
      \[
\lambda^{(+)}_{n} = (2 \pi)^{ -1}  t_{0} ( n-3/4)^{-1} ,\q 
\lambda^{(-)}_{n} = ( 2 \pi)^{ -1}  t_{0} ( n-1/4)^{-1},\q n=1,2,\ldots.
\]
\end{remark}

 \medskip

 {\bf 6.4.}
  Now we are in a position to obtain the asymptotics of the spectrum of Hankel operators whose kernels (or their derivatives) have   jumps of continuity. We combine   Lemma~\ref{powerex2} with the result by V.~V.~Peller (see Theorem~7.4 in Chapter~6 of his book \cite{Pe}) which implies that singular numbers $s_{n} (V)$ of 
    a Hankel operator  $V$ satisfy the bound
    \[
    s_{n}(V)= o(n^{-l-1})
    \]
      if $V$ has a symbol belonging to the Besov class ${\bf B}^{l+1}_{(l+1)^{-1} }({\Bbb R})$. Applying the Weyl theorem on the preservation of the power asymptotics for the sum of operators, we can state the following result.
 
 \begin{theorem}\label{power}
Let $l\in{\Bbb Z}_{+}$, and let $v(t)$ be  the Fourier transform of a function in ${\bf B}^{l+1}_{(l+1)^{-1} }({\Bbb R}) $. Set 
$$h(t)= h_{0} (t_{0}-t)^l+v(t)$$
 for $t\leq t_{0}$ and 
$h(t)=  v(t)$ for $t > t_{0}$. Then eigenvalues of the Hankel operator $H$ have the asymptotics
\[
\lambda_{n}^{(\pm)}=   | h_{0}|  l!  (2\pi)^{-l-1}  t_{0}^{l+1} n^{-l-1} (1+ o(1))
\]
as $n\to \infty$.
 \end{theorem}

 We emphasize that under the assumptions of this theorem the leading terms of the asymptotics of positive and negative eigenvalues are the same.     Of course if $h(t)$ becomes smoother ($l$ increases), then eigenvalues of the Hankel operator $H$ decrease faster as $n\to \infty$.

           %************************************************************
  %************************************************************
\section{Perturbations of the Carleman operator}

%****%****%****%****%****
%****%****%****%****%****%****  

In this section we consider operators $H=H_{0}+V$ where $H_{0}$
 is the Carleman operator ${\bf C}$  (or a more general operator) and the perturbation $V$ belongs to one of the classes introduced in Section~5. Different objects related to the operator $H_{0}$ will be endowed with the index $``0"$,  and objects related to the operator $V$ will be endowed with the index $``v"$.

  \medskip
 
 {\bf 7.1.}
For perturbations     $V$ of  finite rank, we have the following result.
 
  \begin{theorem}\label{FDH}
 Let the sign-function $s_{0} (x)$ of a Hankel operator $ H_{0}$ be    bounded and positive. If $V$ is a Hankel operator  of   finite rank and $H=H_{0}+V$,  then 
  $$N_{-} (H)=N_{-} (V).$$
  In particular, $H \geq 0$ if and only if $V\geq 0$.
 \end{theorem} 
 
 Applying Theorem~\ref{FDH1} to the operator $V$, we get an  explicit formula  for the total number of negative eigenvalues of the operator $H $. 
 
   \begin{corollary}\label{FDH2}
    Let the kernel $v(t)$ of $V$ be given by the formula
  \[
v(t) =\sum_{m=1}^M  P_{m} (t)e^{- \alpha_{m}  t}
\] 
where $P_{m} (t)$ is a polynomial of degree $K_{m}$.
 Define the numbers ${\cal N}_\pm ^{(m)}$   by formula \e{eq:RXm}. 
  Then $N_{-} (H)$      is given by formula  \e{eq:TN}.
 \end{corollary}

 Since for the Carleman operator ${\bf C}$ the sign-function equals $1$, Theorem~\ref{FDH} applies to $H_{0}={\bf C}$.

  The inequality $N_{-} (H)\leq N_{-} (V)$ is of course obvious because $H_{0}\geq 0$. On the contrary, the opposite inequality $N_{-} (H)\geq N_{-} (V)$ looks surprising because the operator $H_{0} $ which may have the continuous spectrum  is much ``stronger" than the operator 
 $V$ of  finite rank. At a heuristic level the equality $N_{-} (H)= N_{-} (V)$ can be explained by the fact that the supports of the sign-functions $s_{0}(x)$ and $s_{v}(x)$ are essentially disjoint. Very loosely speaking, it means  that the operators $H_{0}$ and $V$ ``live in orthogonal subspaces", and hence the positive operator $H_{0}$ does not affect the negative spectrum of $H$.
    The detailed proof of Theorem~\ref{FDH} as well as  that of Theorem~\ref{FDH1}
    will be given in \cite{Y3}.
    
  \medskip
 
 {\bf 7.2.}
  Let ${\bf C } $ be the Carleman operator, and let   $V$ be the Hankel  operator  with kernel 
 \begin{equation}
v(t)= t^k e^{-\alpha t}, \q \alpha> 0,\q k > -1.
\label{eq:E1v}\end{equation}
The operator $V$ is compact, and hence the essential spectrum $\spec_{\rm ess} (H_{\gamma})$ of the operator 
  \begin{equation}
H_{\gamma}= {\bf C } -\gamma V ,\q \gamma\in {\Bbb R},
\label{eq:Cg}\end{equation}
  coincides with the interval $[0,\pi]$.  
   Since the sign-function of the   operator $\bf C$ equals $1$, the sign-function $s_{\gamma}$  of the operator $H_{\gamma}$ equals
 \[
s_{\gamma} (x)= 1 -\gamma s_{v} (x)
\]
where the   function $ s_{v} (x)$ is given by formula \e{eq:bb}.

Let first $k\in (-1,0)$. Observe that $ s_v (x)$ is continuous for $x<\beta=-\ln\alpha$ and $ s_v (x)\to +\infty$ as $x\to \beta-0$ but $s_{v}\in L^1 ({\Bbb R})$ Thus the function $s_{\gamma} (x)\to-\infty$ as $x\to \beta-0$ for all $\gamma>0$, and hence it follows   from Theorem~\ref{HBy} that the   operator $H_{\gamma}$ has an infinite negative spectrum for all $\gamma>0$.

In the case  $k>0$ we use the formula
 \begin{equation}
b(\xi)=   \d(\xi) + b_{v}(\xi)
\label{eq:E1B}\end{equation}
and apply Corollary~\ref{neqx} 
 with $b_{0}(\xi)=  \d (\xi)$ and $b_{\infty}(\xi)=b_{v} (\xi)$.  Since $b_{0} \in C({\Bbb R})'$ and $b_{\infty}  $ has asymptotics \e{eq:E8},  the   operator $H_{\gamma}$ has a negative spectrum for all $\gamma\neq 0$.

Let us summarize the results obtained.

 \begin{proposition}\label{HKC}
 Let $H_{\gamma}={\bf C}-\gamma V$ where $V$ is the Hankel operator    with kernel \e{eq:E1v}. Then:
 
   $1^0$  The operator  $H_{\gamma}$     has an infinite negative spectrum for all $\gamma>0$  if $k\in (-1,0)$.

$2^0$    The operator  $H_{\gamma}$    has negative eigenvalues for all $\gamma\neq 0$  if  $k>0$.
  \end{proposition}

   \medskip

  {\bf 7.3.}
  The  result below   directly follows     from Theorem~\ref{HBy}. 
 
 \begin{proposition}\label{CV}
Suppose that the sign-function   $s_{v}(x) $  of a Hankel operator $V$  is continuous and $s_{v}(x)\to 0$ as $|x| \to\infty$.    Then
 the   operator $H_{\gamma}$  defined by formula \e{eq:Cg} is positive if and only if
\[
\gamma s_{v} (x) \leq 1, \q \forall x\in {\Bbb R}.
\]
If this condition is not satisfied, then   $H_{\gamma}$ has an infinite negative spectrum.
 \end{proposition}

We note that, by Proposition~\ref{extb},  under the assumption of Proposition~\ref{CV} on the sign-function $s_{v} $ the operator $V$ is compact. Of course this assumption is satisfied if $b_{v}  \in L^1 ({\Bbb R})$.

    \begin{example}\label{HKC1}
Let   $v(t)=e^{-t^r}$ where $r<1$. We have seen in subs.~5.3 that its  sign-function $s_{v}(x) =I_r(x)$ where $I_r(x)$ is integral \e{eq:exp}.  Recall that $I_r(x)$  is a nonnegative continuous function of $x\in {\Bbb R}$ and $I_r(x)\to 0$ as $|x|\to\infty$.  Set
 \[
 \nu_r =\max_{x\in {\Bbb R}} I_r(x).
 \] 
Then $H_{\gamma}\geq 0$ if $\gamma\leq \nu_r^{-1}$ and  the operator  $H_{\gamma}$  has infinite negative spectrum for all $\gamma> \nu_r^{-1}$.
  Using explicit formula \e{eq:expZ} it is easy to calculate $\nu_{1/2}=3 \sqrt{6/\pi}e^{-3/2}$.
  \end{example}

  In the case $r>1$ we use formula \e{eq:E1B}. As shown in subs.~5.3,  
  the modulus of the function $ b_{v}(\xi)$   exponentially grows and the periods of its oscillations tend to zero only logarithmically
 as $| \xi |\to\infty$.     Therefore Theorem~\ref{neq} yields the following result.
  
   \begin{proposition}\label{HKC2}
Let   $v(t)=e^{-t^r}$ where $r>1$.   
Then the operator $H_{\gamma} $ has at least one negative eigenvalue for all $\gamma\neq 0$.
  \end{proposition}

\appendix{}

\section{Proof of Lemma~\ref{CHa}}

Set
 \[
F_{\kappa}^{(n)} = \max_{t\in{\Bbb R}_{+}} \big(\la\ln t\ra^\kappa t^n |F^{(n)} (t)|\big)
  \]
  where for shortness we use the notation $\la x\ra= (1+|x|)$.
  
Let us first consider $(\Phi^* F)(\lambda)$ for   $\lambda \in (- 1,1)=:I$.  We have
\[
\sqrt{2\pi} (\Phi^* F) (\lambda)=    \int_{0}^k F(t) e^{i\lambda t}dt+ \int_k^\infty F(t) e^{i\lambda t}dt  ,  \q k = |\lambda|^{-1/2}.
\] 
  The first integral in the right-hand side      is bounded by $F_{0}^{(0)}  |\lambda|^{-1/2}$ which belongs to $L^1 (I)$.
    In the second integral we integrate by parts
  \begin{equation}
  \int_k^\infty F(t) e^{i\lambda t}dt= i \lambda^{-1} F(k) e^{i\lambda k}+i \lambda^{-1}
    \int_k^\infty F'(t) e^{i\lambda t}dt.
      \label{eq:cc1}\end{equation}
    The first term here is bounded by $C  |\lambda|^{-1} \la\ln \lambda\ra^{-\kappa}F_\kappa^{(0)} $ which belongs to $L^1 (I)$ if $\kappa>1$. The second term   is bounded by 
      \[
  | \lambda|^{-1}  \int_k^\infty t^{-1} \la\ln t\ra^{-\kappa}dt
    F_\kappa^{(1)}
    \leq C  | \lambda|^{-1}    \la\ln \lambda\ra^{-\kappa+1} 
    F_\kappa^{(1)} .
 \] 
    It belongs to $L^1 (I)$ if $\kappa>2$.
    
     Next, we consider $(\Phi^* F)(\lambda)$ for   $|\lambda|\geq 1$.
Integrating by parts, we see that
 \begin{equation}
\sqrt{2\pi} (\Phi^* F) (\lambda)=i\lambda^{-1} \int_{0}^k F' (t) e^{i\lambda t}dt +
i\lambda^{-1} \int_k^\infty F' (t) e^{i\lambda t}dt.
    \label{eq:cc3}\end{equation}
The first term here is bounded by
  \[
   | \lambda|^{-1}  \int_0^k t^{-1} \la\ln t\ra^{-\kappa}dt
    F_1^{(\kappa)}
    \leq C  | \lambda|^{-1}    \la\ln \lambda\ra^{-\kappa+1} 
    F_1^{(\kappa)} .
    \]
     It belongs to $L^1 ({\Bbb R}\setminus I)$ if $\kappa>2$. In the second integral in \e{eq:cc3}\ we  once more      integrate by parts, that is, we use formula \e{eq:cc1} with $F(t)$ replaced by $F'(t)$. The function $\lambda^{-2}F'(k)$ is bounded by $ \lambda^{-2}k^{-1}F_{0}^{(1)}=|\lambda|^{-3/2} F_{0}^{(1)}$. For the second term, we use the estimate
     \[
   \big|  \lambda^{-2}
    \int_k^\infty F''(t) e^{i\lambda t}dt    \big| \leq \lambda^{-2}
    \int_k^\infty t^{-2}dt  F_{0}^{(2)}=  \lambda^{-2}  k^{-1} F_{0}^{(2)}=  |\lambda|^{-3/2}  F_{0}^{(2)}.
    \]
    Therefore  the second term in \e{eq:cc3} also belongs to $L^1 ({\Bbb R}\setminus I)$.

\section{The Gaussian kernel}

Here we return to the Hankel operator $H$ with kernel $h(t)=e^{-t^2}$. Now we  proceed from the identity
\begin{equation}
 (Hf,f)= (Q \psi, \psi),
\label{eq:SN1}\end{equation}
 where $\psi(t)=e^{-t^2} f(t)$ and $Q$ is the   integral operator  with real kernel $e^{-2 ts}$ in the space $L^2 ({\Bbb R}_{+})$. We shall use \e{eq:SN1} essentially in the same way as the main identity \e{eq:MAID}. Observe that 
 operators $Q$ with kernels $q(ts)$ which depend only on the product of variables   can be explicitly diagonalized by the Mellin transform (see \cite{Y}). Under fairly general assumptions on
   $q$ the spectrum of $Q$ consists of the interval $[-\gamma ,\gamma ]$ where 
   \[
   \gamma =\sqrt{2\pi}\max_{\xi\in{\Bbb R}}|({\bf M} q)(\xi)|
   \]
   and ${\bf M}$ is the Mellin transform \e{eq:M1}. In particular, for  $q(t)=e^{-2 t}$ the spectrum of $Q$ is absolutely continuous, simple and coincides with the interval  $[-\sqrt{\pi/2 }, \sqrt{\pi/2 } ]$.
     This allows us to check the following assertion.
   
     \begin{proposition}\label{BB}
     The Hankel operator $H$ with kernel $h(t)=e^{-t^2}$ has infinite number of
positive and negative eigenvalues.
  \end{proposition}

    \begin{pf}
Choose some $\mu\in (0,  \sqrt{\pi/2 })$.   For an arbitrary $N$, let $\Delta_{1}^{(+)}, \ldots, \Delta_{N}^{(+)}\subset (  \mu, \sqrt{\pi/2 }  )$ and 
$\Delta_{1}^{(-)}, \ldots, \Delta_{N}^{(-)}\subset ( -\sqrt{\pi/2 },-\mu  )$ be closed mutually disjoint intervals. Choose functions $\varphi_j^{(\pm)}$ such that
  $\varphi_j^{(\pm)}=E_{Q} (\Delta_{j}^{(\pm)})\varphi_j^{(\pm)}$ and $\|\varphi_j^{(\pm)}\|=1$, $j=1,\ldots N$.
 Let $\varphi^{(\pm)}=\sum_{j=1}^N \alpha_{j}\varphi_j^{(\pm)}$ be a linear combination of the functions $\varphi_1^{(\pm)}, \ldots, \varphi_N^{(\pm)}$. Then
 \begin{equation}
 \pm   (Q \varphi^{(\pm)}, \varphi^{(\pm)})= \pm \sum_{j=1}^N |\alpha_{j}|^2  (Q \varphi_{j}^{(\pm)}, \varphi_{j}^{(\pm)})
 \geq \mu \sum_{j=1}^N |\alpha_{j}|^2  \| \varphi_{j}^{(\pm)}\|^2
 = \mu   \| \varphi^{(\pm)}\|^2.
\label{eq:QQ}\end{equation}

 For an arbitrary $\varepsilon>0$, we can choose $\psi_j^{(\pm)}\in C_{0}^\infty ({\Bbb R}_{+})$   such that 
    $\| \psi_j^{(\pm)}-\varphi_j^{(\pm)}\|<\varepsilon$ for all $j=1,\ldots N$. Since the functions $\varphi_j^{(\pm)}$ are orthogonal,  the functions $\psi_j^{(\pm)}$ are linearly independent if $\varepsilon$ is small enough. Moreover, it follows from \e{eq:QQ} that
\[    
 \pm   (Q \psi^{(\pm)}, \psi^{(\pm)})\geq  2^{-1} \mu   \| \psi^{(\pm)}\|^2
\]
 if $\psi^{(\pm)}=\sum_{j=1}^N \alpha_{j} \psi_j^{(\pm)}$ and $\varepsilon$ is small.
 
 Set now  $f ^{(\pm)}(t)=e^{t^2} \psi^{(\pm)}(t)$. Then $f ^{(\pm)}\in L^2 ({\Bbb R}_{+})$ and according to relation  \e{eq:SN1} we have the inequality $\pm (H f_{\pm}, f_{\pm})>0$ on linear subspace of dimension $N$ (except $f_{\pm}=0$).
    Hence the Hankel operator $H$ with kernel \e{eq:F5} has at least $N$ positive and $N$ negative eigenvalues. Since $N$ is arbitrary,  this concludes the proof.
     \end{pf}
    
We emphasize that the operator $H$ is compact 
 while the operator $Q$ has the continuous spectrum. Nevertheless the multiplicities of their positive and negative spectra are the same (infinite).

%%%%%%%%%%%%%%%%%%%%%%%%%%%%%%%%%%%%%%%%
%%%%%%%%%%%%%%%%%%%%%%%%%%%%%%%%%%%%%%%%

 \end{document}